\documentclass[hidelinks,onefignum,onetabnum]{siamart250211}

\usepackage{lipsum}
\usepackage{amsfonts}
\usepackage{graphicx}
\usepackage{epstopdf}
\usepackage{algorithmic}
\usepackage{float}
\ifpdf
  \DeclareGraphicsExtensions{.eps,.pdf,.png,.jpg}
\else
  \DeclareGraphicsExtensions{.eps}
\fi

\newsiamremark{remark}{Remark}
\newsiamremark{hypothesis}{Hypothesis}
\crefname{hypothesis}{Hypothesis}{Hypotheses}
\newsiamthm{claim}{Claim}
\newsiamremark{fact}{Fact}
\crefname{fact}{Fact}{Facts}

\headers{Pseudo-time Continuation}{J. Hart, A. Alexanderian, and B. van Bloemen Waanders}

\title{Preconditioned Pseudo-time Continuation for Parameterized Inverse Problems\thanks{Submitted to the editors DATE.
}}

\author{Joseph Hart\thanks{Department of Scientific Machine Learning, Sandia National Laboratories, Albuquerque, NM 
  (\email{joshart@sandia.gov}, \url{https://www.sandia.gov/ccr/staff/joseph-lee-hart/}).}
\and Alen Alexanderian\thanks{Department of Mathematics, North Carolina State University, Raleigh, NC 
  (\email{aalexan3@ncsu.edu}, \url{https://aalexan3.math.ncsu.edu}).}
\and Bart van Bloemen Waanders\thanks{Department of Scientific Machine Learning, Sandia National Laboratories, Albuquerque, NM 
  (\email{bartv@sandia.gov}, \url{https://www.sandia.gov/ccr/staff/bart-g-van-bloemen-waanders/).}}}

\usepackage{amsopn}

\usepackage{hyperref}
\usepackage{amssymb}
\usepackage{amsmath}
\usepackage{mathtools}

\renewcommand{\vec}[1]{{\mathchoice
                     {\mbox{\boldmath$\displaystyle{#1}$}}
                     {\mbox{\boldmath$\textstyle{#1}$}}
                     {\mbox{\boldmath$\scriptstyle{#1}$}}
                     {\mbox{\boldmath$\scriptscriptstyle{#1}$}}}}

\newcommand{\mat}[1]{\vec{\mathrm{#1}}}

\newcommand{\ta}{\boldsymbol{\theta}}

\newcommand{\tab}{\boldsymbol{\overline{\theta}}}
\newcommand{\tat}{\boldsymbol{\tilde{\theta}}}

\newcommand{\m}{\vec{m}}
\newcommand{\ms}{\vec{m}^*}

\newcommand{\uu}{\vec{u}}
\newcommand{\cc}{\vec{c}}

\newcommand{\II}{\mat{I}}
\newcommand{\HM}{\mat{H}_\text{M}}
\newcommand{\HR}{\mat{R}}

\newcommand{\ipr}[2]{\left\langle {#1}, {#2}\right\rangle_{\HR}}

\newtheorem{example}{Example}[section]
\newtheorem{assumption}{Assumption}[section]

\renewcommand{\H}{\mat{H}}
\newcommand{\B}{\mat{B}}
\newcommand{\E}{\mat{E}}
\newcommand{\F}{\mat{F}}
\newcommand{\W}{\mat{W}}

\newcommand{\R}{\mathbb{R}}

\newcommand{\e}{\boldsymbol{e}}

\newcommand{\f}{\boldsymbol{f}}

\ifpdf
\hypersetup{
  pdftitle={Preconditioned Pseudo-time Continuation for Parameterized Inverse Problems},
  pdfauthor={J. Hart, A. Alexanderian, and B. van Bloemen Waanders}
}
\fi

\begin{document}

\maketitle

\begin{abstract}
We consider parameterized variational inverse problems that are constrained by
partial differential equations (PDEs). We seek to efficiently compute the
solution of the inverse problem when auxiliary model parameters, which appear in
the governing PDE, are varied. Computing the solution of the inverse problem for
different auxiliary parameter values is crucial for uncertainty quantification.
This, however, is computationally challenging since it requires
solving many optimization problems for different realizations of the 
auxiliary parameters. We
leverage pseudo-time continuation and
solve an initial value problem to evolve the optimal solution along an auxiliary parameter path. This
article introduces the use of an adaptive quasi-Newton Hessian preconditioner to
accelerate the computation. Our proposed preconditioner exploits properties of
the pseudo-time continuation process to achieve reliable and efficient
computation. We elaborate our proposed framework and elucidate its properties
for two nonlinear inverse problems. 

\end{abstract}



\section{Introduction}

Inverse problems arise in a variety of applications where a quantity of interest
cannot be measured directly. Rather, indirect observations 
are made and the quantity of interest is related to the observations via a
mathematical model. Often, the model is a partial differential equation (PDE),
the inversion parameter of interest is a model parameter, and the observable
quantity is a function of the solution of the PDE.
Various approaches exist to estimate the model parameters.
All the approaches require solving the PDE for different model parameters and
comparing the model output with the observed data by evaluating a data-misfit functional.

The PDE governing an inverse problem typically has a number of \emph{auxiliary}
parameters.  These are  model parameters that are needed to fully specify the
governing PDE, but are not the inversion parameters.
We focus on the challenge of understanding the 
impact of auxiliary parameters on the inverse problem.
This is crucial for quantifying
the uncertainty in the estimated inversion parameter. A first-principles approach is to
jointly estimate the inversion parameter of interest and the auxiliary
parameters~\cite{Crestel_2019}. However, this is frequently intractable due to
the increased dimensionality induced by including the auxiliary parameters in the
inversion. The auxiliary parameter dimension may be decreased via
screening and sensitivity analysis techniques~\cite{future_of_sa,gsa_primer}. 
Although such approaches may provide auxiliary parameter dimension reduction,
the presence of influential auxiliary parameters is frequently inescapable and
thus must be accounted for in the inverse problem. 

In a Bayesian paradigm, auxiliary parameter uncertainty can be accounted for by
the Bayesian approximation error (BAE)
approach~\cite{Kaipio_2008,bayes_approx_error_petra}. This is achieved by
pre-marginalizing the auxiliary parameter uncertainty, which results in a
modified noise model.  In this article, we focus on variational inverse
problems.  Such deterministic inverse problems require solving an optimization
problem whose objective is the data-misfit---the difference between the model
prediction and the observed data. Also, in ill-posed inverse problems, which are
our primary focus, a regularization term is added to the data-misfit functional
to define the objective function. 

A pragmatic approach is to fix the auxiliary parameter to a nominal
value and solve the optimization problem to determine the parameter of interest
that results in the model agreeing with the observational data. However, to
incorporate auxiliary parameter uncertainty, it is necessary to solve many such
optimization problems for different auxiliary parameters. This is
computationally prohibitive due to the many optimization solves required. 

To avoid solving optimization problems repeatedly, hyper-differential sensitivity analysis
(HDSA) was developed to compute the sensitivity of the nominal optimization
problem's solution with respect to perturbations of the auxiliary
parameters~\cite{hart_2021_bayes,HartvanBloemenWaanders20,sunseri_2,sunseri_hdsa}. This provides a
linear approximation of the auxiliary parameter to optimal solution mapping. 
The key benefit of HDSA is its computational scalability.
Namely, HDSA provides an efficient and accurate framework for 
estimating the sensitivity of
the estimated inversion parameter with respect to the auxiliary parameter
in large-scale inverse problems with high-dimensional inversion and 
auxiliary parameters.
However, this linear approximation may be insufficient for large
perturbations, e.g., when there is significant auxiliary parameter uncertainty.

To overcome the limitations of linear approximation, continuation and
path-following algorithms can be used to solve perturbed optimization problems.
Early work focused on continuation methods to solve systems of nonlinear
equations. See the books~\cite{allgower_georg,allgower_georg_siam}, and
the references therein, for foundational contributions in continuation methods,
such as the predictor-corrector algorithm. These ideas also translate to optimization
problems since the first-order optimality condition defines a system of
nonlinear equations. Furthermore, optimization problems have additional
structure, such as symmetry of the Hessian matrix, that can be exploited to
develop efficient algorithms. See the book~\cite{param_opt_book} for an
exposition of continuation methods for parametric optimization problems
and~\cite{Kungurtsev_2017} for a treatment of dual degenerate problems. 
See also the recent work~\cite{Liu_2025,Si_2024}, which demonstrates the utility of
path-following---continuation along a path through the parameter space---to
compute optimal parameters for statistical models. 
These works focus on
regression problems where the primary source of complexity is in the
regularization term arising from a Bayesian prior. 
On the other hand, the present work
focuses on inverse problems where the greatest challenge arises from the
PDE-constrained data-misfit term.

We build on the recent work~\cite{AlexanderianHartStevens23} that introduced the
idea of pseudo-time continuation, which is synonymous with path-following, in
the context of parameterized PDE-constrained inverse problems.  The present
article differs from the breadth of continuation literature as a result of the
computational challenges and mathematical structures that arise from
PDE-constrained optimization. These include the necessity of iterative solvers
for Hessian systems and the computation of Hessian-vector products via the
adjoint method. Our contributions in this article are: (i) the design of
adaptive Hessian preconditioners that improve computational efficiency, and
(ii) demonstration on challenging PDE-constrained inverse problems that
elucidate the key features of the proposed approach. We focus on efficiently
computing the minimizer for one specific auxiliary parameter perturbation
by tracking minimizers along a path from the nominal parameter value to the
perturbed parameter value. There are natural extensions of this work to
computing minimizers for many different parameter perturbations. This is 
discussed in our conclusions. 

The article is organized as follows. Section~\ref{sec:prelim} provides the
mathematical framing of the problem and recalls relevant literature. The
computational kernel of pseudo-time continuation is iterative linear solves to
invert the Hessian. Our approach to preconditioning these solves is introduced
in Section~\ref{sec:preconditioners}. Numerical results are presented in
Section~\ref{sec:numerical_results} to demonstrate the benefits of our proposed
approach and give insight into the effect of the algorithm parameters (step
sizes, tolerances, etc.). In Section~\ref{sec:conc}, we provide concluding remarks 
and discuss potential extensions 
of the proposed approach.  

\section{Preliminaries}\label{sec:prelim}
Consider parameterized variational inverse 
problems of the form
\begin{equation}\label{equ:optim}
    \min_{\m \in \R^n} J(\m, \ta) \coloneq \frac12 \| \mat{F}(\m; \ta) - \vec{y} \|^2 + \frac12 \m^\top \mat{R} \m. 
\end{equation}
Here, $\m$ is the vector of inversion parameters, $\mat{F}$ is a
parameter-to-observable map, $\vec{y} \in \R^d$ is measurement data, $\ta$
is a vector of auxiliary parameters, and $\mat{R} \in \R^{n \times n}$ is a symmetric positive definite matrix that is used for regularization.  We consider inverse problems governed by
PDEs in which the auxiliary parameters correspond to model parameters that are
not being estimated. This is why the parameter-to-observable map is written as a
function of both $\m$ and $\ta$.  The governing PDE must be solved each time
$\mat{F}(\m; \ta)$ is evaluated.  Consequently, solving~\eqref{equ:optim} is
challenging due to the computational cost of solving the PDE repeatedly as
$\mat{F}(\m; \ta)$ is evaluated for different inputs. 

We assume throughout that $J(\m,\ta)$ is twice continuously differentiable with respect to $(\m,\ta)$
and use $\nabla_\m J$ and $\nabla_{\m,\m} J$ to denote the gradient and Hessian of $J$ with respect to $\m$; 
similar notation is used to denote the $\ta$ derivatives. We consider problems for which matrix-vector products with the Hessian $\nabla_{\m,\m} J(\m,\ta)$ may be computed and use Newton methods to solve~\eqref{equ:optim}. The Hessian information
facilitates rapid convergence of the optimization iterates. Given a fixed $\ta$ 
and current optimization iterate $\m_k$, the Newton search direction is obtained by solving
the system of linear equations
\begin{align} \label{eqn:newton_system}
\H(\m_k,\ta) \vec{s}_k = - \nabla_\m J(\m_k,\ta),
\end{align}
where $\H(\m_k,\ta) = \nabla_{\m,\m} J(\m_k,\ta)$ is the Hessian
of the objective function. The Newton step is $\m_{k+1}=\m_k + \vec{s}_k$. A
globalization strategy, such as line search or trust region, is necessary to
ensure convergence of the optimization algorithm. However, when sufficiently
close to the minimizer, the Newton steps result in rapid (quadratic) convergence
of the optimization iterations toward the minimizer. 

When solving the inverse problem, traditionally $\ta$ is fixed at a vector
$\tab$ of nominal auxiliary parameters. When the auxiliary parameters are subject 
to considerable uncertainty, it is important to quantify the uncertainty 
in the solution of the inverse problem. However, doing this by repeatedly 
solving~\eqref{equ:optim} for different values of $\ta$ is computationally expensive. 
In what follows, we let $\ms(\ta)$ denote a minimizer of $J$. Note that for non-convex optimization problems,
uniqueness of $\ms(\ta)$ is ensured by restricting to a neighborhood of the local minimizer computed
when $\ta=\tab$; see~\cite{AlexanderianHartStevens23} for details. 
In this section, we recall the pseudo-time continuation approach and the predictor-corrector algorithm 
to compute $\ms(\ta)$ with $\ms(\tab)$ used as a starting point. Section~\ref{sec:preconditioners} details this article's contribution
to enhance the algorithms through judicious preconditioning to compute $\ms(\ta)$ more efficiently.

\subsection{Pseudo-time continuation}
Consider the Jacobian of $\ms$ with respect to $\ta$, which we denote by $\nabla_{\ta} \ms$.
Applying the Implicit Function Theorem to the first-order optimality condition for~\eqref{equ:optim} yields 
\begin{equation*}
		\nabla_{\ta} \ms(\ta) = -\H(\ms(\ta),\ta)^{-1}
                                 \B(\ms(\ta),\ta).
\end{equation*}
Here $\B(\ms(\ta),\ta) = \nabla_{\m,\ta} J(\ms(\ta),\ta)$ is the mixed second derivative of the objective function and $\H(\ms(\ta),\ta)$ is the Hessian previously introduced in the discussion on Newton's method; see~\eqref{eqn:newton_system}. 
We can consider a linear approximation model 
\begin{equation} \label{equ:linear_approx}
\ms(\ta) \approx \ms(\tab) + \nabla_{\ta} \ms(\tab)(\ta - \tab).
\end{equation}

This approximation is only valid in a
neighborhood of $\tab$ and may therefore be inaccurate if $\| \ta - \tab \|$ is large. To
overcome this limitation, we consider a pseudo-time continuation approach. Let
$\tat$ be a parameter vector for which we seek to compute $\ms(\tat)$ and
parameterize a path from $\tab$ to $\tat$ via $\ta(t) = \tab + (\tat-\tab)t$, $t
\in [0,1]$. With an abuse of notation, we consider $\ms(t) \coloneqq
\ms(\ta(t))$. It follows that 
\[
\frac{d \ms}{dt}(t) = \frac{\partial \ms}{\partial \ta}(\ta(t)) \frac{d \ta}{dt}(t) = -\H(\ms(\ta(t)),\ta(t)) ^{-1} \B(\ms(\ta(t)),\ta(t)) (\tat-\tab).
\]
Let 
\begin{align} \label{eqn:f}
\f(\ms(t),\ta(t)) \coloneq -\H(\ms(\ta(t)),\ta(t)) ^{-1} \B(\ms(\ta(t)),\ta(t)) \Delta \ta,
\end{align}
where $\Delta \ta \coloneq \tat - \tab$. We may 
compute $\ms(\tat)$ by solving the initial value problem 
\begin{align}
\label{equ:ms_ode}
& \frac{d\ms}{dt} = \vec{f}(\ms(t),\ta(t)) \quad t \in [0, 1], \\
& \ms(0) = \ms(\tab). \nonumber
\end{align}

To ensure that $\H(\ms(\ta(t)),\ta(t))$ is invertible, and that $\ms(\ta(t))$ is tracking a strict minimizer, we make the following assumption:
\begin{assumption}
$\H(\ms(\ta(t)),\ta(t))$ is positive definite for all $t \in [0,1]$.
\end{assumption}
In general, tracking a stationary point $\ms(\ta(t))$ that satisfies $\nabla_\m J(\ms(\ta(t)),\ta(t))=\vec{0}$ does not guarantee positive definiteness of the Hessian. However, this assumption is reasonable in practice for many inverse problems arising from physical systems, particularly since regularization is typically used in ill-posed problems. Lastly, we note that a common approach to addressing near-zero eigenvalues in the Hessian is to project onto a likelihood informed subspace~\cite{cui_2014}. If indefiniteness is expected at the nominal minimizer $\ms(0)$, we can consider optimization in only the likelihood informed subspace, which will ensure satisfaction of the assumption in almost all practical scenarios.

Observe that
the linear approximation~\eqref{equ:linear_approx} corresponds to discretizing~\eqref{equ:ms_ode} 
with forward Euler and integrating with a single time step of size 1. However, by viewing the evolution of the
minimizer via an ODE system~\eqref{equ:ms_ode}, we may consider higher-fidelity time stepping
algorithms to ensure more accurate solutions.
To demonstrate the benefits of pseudo-time
continuation relative to an optimization algorithm, we consider a simple example below.
\begin{example} \label{example:illustrative}
Consider the parameterized optimization problem
\begin{align} \label{equ:illustrative_example}
    \min_{m \in \R} J(m, \theta) \coloneq (m-\theta)^6 + 0.01m^2
\end{align}
where $\theta \in \R$ is an uncertain auxiliary parameter. This illustrative example captures
the common inverse problem scenario that the objective function has small curvature in
some parameter directions~\cite{GhattasWilcox21}. We assume
that~\eqref{equ:illustrative_example} has been solved when $\theta=1.0$ to
produce the minimizer $m^* = 0.702$. To compute the minimizer when $\theta=4.0$, we
take $m=0.702$ as the initial iterate and use Newton's method to
solve~\eqref{equ:illustrative_example} with $\theta=4.0$. For comparison, we 
use pseudo-time continuation taking a forward Euler
discretization of~\eqref{equ:ms_ode} with $3$ time steps. Figure~\ref{fig:illustrative_example_post_opt_history} displays the 
iterations produced by Newton's method and the time steps produced by pseudo-time continuation. 

\begin{figure}[ht]
\centering
    \includegraphics[width=0.49\textwidth]{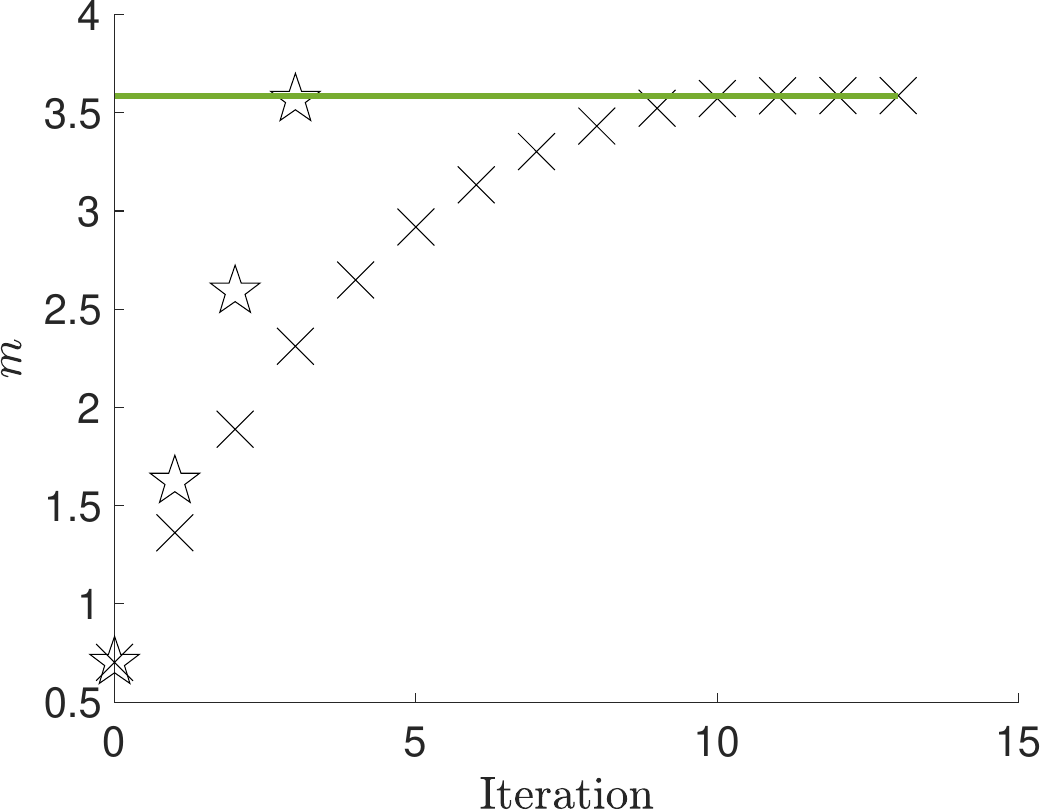}
    \caption{Newton's method iterations and pseudo-time continuation steps generated by minimizing $J(m,4.0)$ when initialized with the minimizer of $J(m,1.0)$. The iterations for Newton's method and time steps for pseudo-time continuation approach are marked by crosses and stars, respectively. The horizontal line indicates the minimizer of $J(m,4.0)$.}
  \label{fig:illustrative_example_post_opt_history}
\end{figure}

We note that the pseudo-time continuation approach converges faster than Newton's method. We
conjecture that this improved convergence occurs for the following interrelated reasons:
\begin{enumerate}

\item 
In Newton's method, we consider a sequence of iterations that minimize
successive local quadratic approximations of the objective function.  On the
other hand, in pseudo-time continuation we consider \eqref{equ:ms_ode}, which
describes the evolution of the minimizer as a function of $\theta$.  Subsequently,
we approximate the numerical solution of that initial value problem with a
time-stepping method.  Hence, the pseudo-time continuation approach permutes the
actions of approximation and minimization. 

\item By taking $3$ forward Euler time steps, we perturb $\theta$ by $1.0$ at
each step whereas Newton's method takes a size $3.0$ perturbation ($1.0 \to
4.0$) of $\theta$ initially. 

\item Since our pseudo-time continuation algorithms changes $\theta$ gradually,
we stay within the Newton basin such that a time step of the
form~\eqref{eqn:fe_time_step} tracks the minimizer well.  Whereas initializing
$J(m,4.0)$ at $m=0.702$ is outside the Newton basin and hence does not enjoy
the rapid quadratic convergence that characterizes Newton's method when it is
``close enough" to the minimizer. 
\end{enumerate}
\end{example}

These observations motivate the study of pseudo-time continuation algorithms.
Also, as demonstrated shortly, posing the perturbed optimization problem in
terms of the initial value problem~\eqref{equ:ms_ode} provides other
opportunities to accelerate computation.

\subsection{Predictor-corrector algorithm}
To obtain the minimizer corresponding to a perturbed 
auxiliary parameter vector $\tat$, we can, in principle, 
solve the initial value problem~\eqref{equ:ms_ode}
using a time-stepping method. However, the ODE system~\eqref{equ:ms_ode} has additional structure
that is not exploited by time-stepping methods. Specifically, for any $t$,
$\ms(t)$ is a local minimizer of $J(\m,\ta(t))$ and hence $\nabla_\m J(\vec{m}(t),\ta(t)) = \vec{0}$. If we integrate~\eqref{equ:ms_ode}
via a traditional time-stepping method, the accumulation of error results in $\| \nabla_\m J(\vec{m}(t),\ta(t))  \|$ increasing as
a function of time. Rather, we can use a predictor-corrector procedure~\cite{allgower_georg,allgower_georg_siam,param_opt_book} where 
a traditional time-stepping method predicts the solution and a Newton optimization iteration corrects the predicted solution to ensure that 
$\| \nabla_\m J(\vec{m}(t),\ta(t))  \|$ is less than a specified tolerance. Note that despite using the
terminology, and having many conceptual similarities, predictor-corrector algorithms in the present
context are different from the classical predictor-corrector ODE solvers. In particular, predictor-corrector ODE solvers
use two different time stepping schemes, whereas the algorithms we consider use a time stepping scheme for prediction followed by a Newton step for correction. This corrector is unique to path-following of optimization solutions.

To facilitate the discussion, we first consider 
the setup of a time-stepping approach for~\eqref{equ:ms_ode}.
We focus on explicit methods and discretize~\eqref{equ:ms_ode} using $N$ time steps.
Let $\Delta t = 1/N$ denote the step size 
and let $t_k = k \Delta t$, for $k \in \{0, \ldots, N\}$.
We let $\m_k$ denote the estimated
value of $\ms(t_k)$ and 
$\ta_k=\ta(t_k)$ for $k \in \{0,1,\ldots,N\}$. 

Consider the numerical solution of~\eqref{equ:ms_ode}
approximated using
forward Euler,
\begin{equation} \label{eqn:fe_time_step}
\m_{k+1} = \m_k - \Delta t \H(\m_k,\ta_k)^{-1} \B(\m_k,\ta_k)\Delta\ta , \quad k = 0, \ldots, N-1.
\end{equation}

The Euler-Newton predictor-corrector scheme consists of the steps:
\begin{subequations}\label{equ:predictor_corrector}
	\begin{align}
	\m_{k+1}^\text{pred} &= \m_k - \Delta t \H(\m_k,\ta_k)^{-1} \B(\m_k,\ta_k)\Delta\ta
	\label{equ:predictor_main}
	\\
	\m_{k+1}     &= \m_{k+1}^\text{pred} - \H(\m_{k+1}^\text{pred},\ta_{k+1})^{-1} \nabla_\m J(\m_{k+1}^\text{pred},\ta_{k+1}). 
	\label{equ:corrector_main}
	\end{align}
\end{subequations}
In the above procedure, the forward Euler 
step~\eqref{equ:predictor_main} serves as a predictor. Subsequently, 
one step of Newton's method is used as a corrector. 
In practice, the predictor step results in an approximation to 
$\m_{k+1}$ that falls close to the basin of attraction for the 
$\ms(t_{k+1})$. The corrector step then uses this as an initial guess 
for a Newton step.
As noted later in the article, additional corrector steps may be needed to
ensure optimality within a prescribed tolerance.  However, the
optimality typically holds with one or two corrector steps thanks to the rapid
convergence of Newton's method. 

The Euler-Newton predictor-corrector scheme~\eqref{equ:predictor_corrector} is commonly used~\cite{allgower_georg_siam,param_opt_book}. To provide further intuition, we recall classical observations of similarity between the forward Euler predictor step~\eqref{equ:predictor_main} and an optimization iterate using Newton's method. Specifically, a first-order approximation of $\Delta t \B(\m_k,\ta_k)\Delta\ta$ yields 
\begin{align*}
\Delta t \B(\m_k,\ta_k)\Delta\ta & \approx
\nabla_\m J(\m_k,\ta_k+\Delta t \Delta \ta) - \nabla_\m J(\m_k,\ta_k) \\
& \approx  \nabla_\m J(\m_k,\ta_k+\Delta t \Delta \ta),
\end{align*}
where the latter approximation follows since $\nabla_\m J(\m(t_k),\ta_k)=\vec{0}$ and $\m_k \approx \ms(t_k)$. Consequently, the forward Euler step~\eqref{equ:predictor_main} is approximately given by
\begin{align*}
 \m_k - \H(\m_k,\ta_k)^{-1} \nabla_\m J(\m_k,\ta_{k+1}).
\end{align*}
To minimize $J(\m,\ta_{k+1})$ with Newton's method using an initial iterate $\m_k$, the iteration takes the form
\begin{align*}
 \m_k - \H(\m_k,\ta_{k+1})^{-1} \nabla_\m J(\m_k,\ta_{k+1}).
\end{align*}
Hence, the predictor step may be viewed as a type of Newton step that is informed by the $\ta$ variability.

There is a lot to say about the procedure outlined 
in~\eqref{equ:predictor_corrector}.
In the first place, one may consider higher-order
time-stepping methods for the predictor.
In Section~\ref{sec:predictor_step}, 
we present and discuss the advantages of a second-order accurate method 
in the predictor step. Another important consideration is the cost of 
the Hessian solves in the predictor and corrector steps. In the 
case of optimization problems arising from discretized inverse problems 
constrained by PDEs, these Hessian solves are expensive and are 
tackled using Krylov iterative methods. Specifically, we rely on 
the Conjugate Gradient (CG) method. The cost of 
CG solves for these Hessian systems is dominated by expensive 
matrix-vector products with the 
Hessian---each Hessian-vector product requires two PDE solves. 
In this context, preconditioning is an essential tool 
to reduce the number of CG iterations. We outline our 
proposed preconditioners in Section~\ref{sec:preconditioners}. Our approach leverages properties of the pseudo-time continuation
and predictor-corrector algorithms to evolve the preconditioner with time.

\subsection{The predictor step}\label{sec:predictor_step}

The role of the predictor step is 
to provide a sufficiently accurate initial guess for the corrector step 
to ensure rapid convergence. 
Much of the work on optimal solution path-following predictor-corrector methods rely on
forward Euler time-stepping~\cite{allgower_georg_siam}. However, if the auxiliary parameter perturbation 
is large, we may need a small time step 
when using forward Euler.

Forward Euler is first-order accurate, that is,
the truncation error of time discretization is $\mathcal O(\Delta t)$. A second-order method has a truncation
error that is $\mathcal O(\Delta t^2)$. This has the potential to reduce computational cost by taking larger time steps.
There are multiple second-order methods that we can
consider. In this article, we focus on a modified Euler method defined by the
steps:
\begin{subequations}\label{equ:mod_euler}
	\begin{align}
& \m_{k+\frac{1}{2}} = \m_k + \frac12 \Delta t \f(\m_k,\ta_k) \label{eqn:mod_euler_1} \\
& \m_{k+1} = \m_k +  \Delta t \f(\m_{k+\frac{1}{2}},\ta_{k+\frac{1}{2}}). \label{eqn:mod_euler_2}
	\end{align}
\end{subequations}

This modified Euler method corresponds to taking a step of size $\frac 12 \Delta t$
using forward Euler~\eqref{eqn:mod_euler_1} to predict an intermediate
$\m_{k+\frac{1}{2}}$, followed by a step of size $\Delta t$ from time $t_k$ to
$t_{k+1}$ using the slope $\f(\m_{k+\frac{1}{2}},\ta_{k+\frac{1}{2}})$ evaluated
at the midpoint predicted by the forward Euler step~\eqref{eqn:mod_euler_2}. There are many potential
second-order numerical integration schemes, see~\cite{Gautschi_book} for an overview.
This modified Euler method is preferred to other second-order schemes because it increments time by $\frac{1}{2} \Delta t$ in both
~\eqref{eqn:mod_euler_1} and~\eqref{eqn:mod_euler_2}.

Consider forward Euler time-stepping with $N$ time
steps and modified Euler time-stepping with $M=\frac{N}{2}$ time steps (assuming $N$ is even).
The time nodes
$t_0,t_\frac{1}{2},t_1,t_{1+\frac{1}{2}},\dots,t_{M-\frac{1}{2}},t_{M}$
in the modified Euler time discretization correspond exactly to the time nodes
$t_0,t_1,\dots,t_N$ in the forward Euler time discretization. Hence, both
methods require $N$ evaluations of $\f$ and thus have approximately the same
computational cost. However, the truncation error of forward Euler is $\mathcal
O(\frac{1}{N})$ whereas the truncation error of modified Euler is $\mathcal
O(\frac{1}{M^2})=\mathcal O(\frac{4}{N^2})$. If $N>4$ and the constants hidden in the $\mathcal O$ are equal, then the modified Euler method has a smaller truncation error with the same computational cost.
This simple analysis 
is instructive to see  when and
why higher-order methods may be superior. Our numerical results presented in
Section~\ref{sec:numerical_results} demonstrate the improvements that
can be attained using modified Euler when the auxiliary parameter perturbation is large enough to mandate small forward Euler time steps.  Although we omit discussion of
higher-order (greater than second-order) methods now, we highlight
their potential benefits in our concluding remarks; see Section~\ref{sec:conc}.

\section{Preconditioners for the Hessian solves}
\label{sec:preconditioners}

The pseudo-time continuation algorithm requires many evaluations of $\f(\m,\ta)$, which was defined in~\eqref{eqn:f}.  The key
challenge when evaluating $\f$ is computing the action of $\H(\m,\ta) ^{-1}$ on
vectors.  These Hessian solves are performed using the Conjugate Gradient (CG)
method.  
Each CG iteration requires a Hessian apply, which involves solving the 
incremental state and incremental adjoint equations; see Appendix~\ref{appendix:computation} for details. Hence, to achieve
computational efficiency, it is crucial to precondition the Hessian solves so
that we reduce the number of CG iterations.

In this section, we present novel quasi-Newton methods to build effective
Hessian preconditioners within the pseudo-time continuation algorithm. We
note at the onset that this is made possible because our approach starts at a
local minimizer and tracks its changes as $\ta$ is perturbed. Furthermore, as we previously assumed (and argued is reasonable in practice), the Hessian matrices remain positive definite 
throughout the pseudo-time steps. In a general optimization
setting, effective preconditioning of the Hessian system is more difficult. We note
that the idea to design Hessian preconditioners to accelerate CG solves within a path-following algorithm
was also considered in~\cite{Si_2024}. In that work, the optimization problem arose from 
Bayesian hierarchical statistical modeling where the complexity stemmed from the 
prior model. Specialized preconditioners were developed in~\cite{Si_2024} by exploiting
structure associated with the statistical modeling problem. In our context, the complexity 
comes from the underlying PDE constraint. This 
motivates a different approach to design preconditioners that
adapt using information from the PDE (and associated adjoint) solves.

In preparation for our contributions, we first recall the Preconditioned Conjugate
Gradient (PCG) method in Algorithm~\ref{alg:pcg}, and refer the reader
to~\cite{kelley_book_1995} for a more detailed discussion. Note that we present PCG
using Euclidean inner products. In general, if the linear system arises from the
discretization of an linear operator defined on function spaces, then the PCG solver
should use a weighted inner product that arises from discretizing the function space 
inner product~\cite{pcg_prec_book}. However, in the context of the Hessian systems in Newton
optimization algorithms, the inner product weighting is already embedded in the discretized
Hessian so that the algorithms may use Euclidean inner products.

\begin{algorithm} 
\caption{Preconditioned Conjugate Gradient (PCG)}
\begin{algorithmic}[1]
\STATE \textbf{Input: } $\vec{b} \in \R^n$, $\mathbf{A} \in \R^{n \times n}$, $\mathbf{E} \in \R^{n \times n}$, $\epsilon_{\text{CG}}>0$
\STATE \textbf{Initialize: } $\vec{r}=\vec{b}$, $\vec{x}=\vec{0}$, $\rho_0=\| \vec{r} \|^2$, $k=1$
\WHILE{$\sqrt{\rho_{k-1}} > \epsilon_{\text{CG}} \| \vec{b} \|$}
\STATE $\vec{z} = \mathbf{E} \vec{r}$
\STATE $\tau_{k-1} = \vec{z}^\top \vec{r}$
\IF{$k=1$}
\STATE $\beta=0$ and $\vec{p}=\vec{z}$
\ELSE
\STATE $\beta = \tau_{k-1}/ \tau_{k-2}$ and $\vec{p}=\vec{z} + \beta \vec{p}$
\ENDIF
\STATE $\vec{w} = \vec{A} \vec{p}$
\STATE $\alpha = \tau_{k-1}/(\vec{p}^\top \vec{w})$
\STATE $\vec{x} = \vec{x} + \alpha \vec{p}$
\STATE $\vec{r} = \vec{r} - \alpha \vec{w}$
\STATE $\rho_k = \| \vec{r} \|^2$
\STATE $k=k+1$
\ENDWHILE
\end{algorithmic}
\label{alg:pcg}
\end{algorithm}

Let $\mathbf{E}$ denote the preconditioner, which must be symmetric positive definite and should
approximate the inverse Hessian $\H(\m,\ta)^{-1}$. Ensuring that the condition
number of $\mathbf{E} \H(\m,\ta)$ is small improves the convergence of the PCG
iteration. An alternative perspective on preconditioning is that we desire to
have the eigenvalues of $\mathbf{E}^\frac{1}{2}
\H(\m,\ta)\mathbf{E}^\frac{1}{2}$ clustered. In exact arithmetic, if
$\mathbf{E}^\frac{1}{2} \H(\m,\ta) \mathbf{E}^\frac{1}{2}$ has $m$ distinct
eigenvalues then PCG will converge in at most $m$ iterations. From the
perspective of clustering eigenvalues, an attractive preconditioner is
$\mathbf{E}=\mathbf{R}^{-1}$, where $\mathbf{R}$ is the
regularization operator in~\eqref{equ:optim}~\cite{GhattasWilcox21}. To see this, let $\HM(\m,\ta)$ be the Hessian of
the data-misfit term $\frac{1}{2} \| \vec{F}(\m,\ta) - \vec{y} \|^2$ and observe that 
\begin{align*}
\mathbf{R}^{\frac{1}{2}} \H(\m,\ta)\mathbf{R}^{\frac{1}{2}}  = \mathbf{R}^{\frac{1}{2}} \HM(\m,\ta)\mathbf{R}^{\frac{1}{2}}  + \vec{I},
\end{align*}
where $\mathbf{I}$ denotes the identity matrix. In many applications,
$\mathbf{R}^{\frac{1}{2}} \HM(\m,\ta)\mathbf{R}^{\frac{1}{2}}$ is approximately low-rank as a result of
data sparsity and smoothing properties of the parameter-to-observable map~\cite{GhattasWilcox21}. Hence
the regularization-preconditioned data-misfit Hessian $\mathbf{R}^{\frac{1}{2}}
\H(\m,\ta)\mathbf{R}^{\frac{1}{2}} + \mathbf{I}$ will have many of its eigenvalues
clustered near $1$.  
While $\mathbf{R}$ is a commonly used preconditioner 
in practice, we demonstrate that superior convergence can
be achieved by adapting the preconditioner within the pseudo-time continuation
algorithm. To this end, we introduce two
novel quasi-Newton updates that exploit the pseudo-time continuation structure.

Quasi-Newton methods are commonly used for numerical optimization in problems
where gradients are available, but Hessians are not. The idea of a quasi-Newton
method is to begin the optimization algorithm with an initial estimate of the
Hessian (frequently a poor one) and update it at each optimization  iteration using
gradient evaluations. We build on this idea and propose an approach wherein we initialize an
estimate of the Hessian and update it after each time step within the predictor 
corrector iterations~\eqref{equ:predictor_corrector}.

In the subsections that follow, we present two quasi-Newton preconditioner
updates. The first one is based on introducing a secant equation that captures
the $\ta$ variations. We refer to this as the \textit{parametric quasi-Newton
update}. The second one is an adaptation of the block BFGS
method~\cite{gao_2018}, which leverages the PCG iteration history to enhance the
preconditioner in future time steps.  We refer to this second update as the
\textit{block quasi-Newton update}.  Central to both updates is that 
only existing gradients and Hessian-vector products are used to improve the
approximation. Hence, additional PDE solves are not required. However,
storage of additional vectors is needed that would otherwise be removed
from memory.

\subsection{Parametric Quasi-Newton Update}

When solving~\eqref{equ:ms_ode}, we require a sequence of $\H(\m_k,\ta_k)^{-1}$
applies, where $(\m_k,\ta_k)$ changes gradually.  Traditional quasi-Newton
methods, which rely on gradient differences, can be used to precondition these
linear solves.  However, this ignores the dependence of $\m_k$ on $\ta_k$.  In
the traditional setting, the Hessian only depends on $\m$. Here, we derive a
quasi-Newton update that accounts for the parameter variations.

To initialize the preconditioned pseudo-time continuation algorithm, we require
that a preconditioner for $\H(\m_0,\ta_0)$ is provided. 
In the context of inverse
problems, there are two natural options to define an initial preconditioner: (i) use the inverse
of the matrix $\mathbf{R}$ in~\eqref{equ:optim}, 
or (ii) use a low-rank approximation of 
$\vec{R}^{\frac{1}{2}} \HM(\m_0,\ta_0)\vec{R}^{\frac{1}{2}}$ to approximate  
$\H(\m_0,\ta_0)^{-1}=( \vec{R}^{\frac{1}{2}}
\HM(\m_0,\ta_0)\vec{R}^{\frac{1}{2}}+\mathbf{I})^{-1}$. We
provide a detailed derivation of the latter 
in Appendix~\ref{appendix:init_hess}. 
For the analysis that follows, we let $\E_k$ denote our preconditioner that 
approximates $\H(\m_k,\ta_k)^{-1}$, for $k=0,1, \ldots, N$, where
$\E_0$ is provided via one of the options discussed above.

We assume that a preconditioner $\E_{k-1}$ is available and seek to update it to produce a preconditioner $\E_k$ for the next time step.
To derive a quasi-Newton update, we introduce a secant equation in the context
of pseudo-time continuation. Let
\[
\phi_k(\m)=J(\m_k,\ta_k) + (\m-\m_k)^\top \nabla_\m J(\m_k,\ta_k) + \frac{1}{2} (\m-\m_k)^\top \widehat{\vec{H}}_k (\m-\m_k)
\]
be a local quadratic approximation of the objective function $J(\m,\ta_k)$
centered at $\m_k$. Here, $\widehat{\vec{H}}_k$ denotes an approximate Hessian. We 
seek to determine $\widehat{\vec{H}}_k$ and define $\E_k = \widehat{\vec{H}}_k^{-1}$.

In traditional optimization, $\ta_k$ is fixed and only $\m_k$ changes with each iteration, whereas in our context both $\m_k$ and $\ta_k$ change with each time step. We enforce the secant conditions 
\begin{align*}
\nabla_\m \phi_k(\m_k)=\nabla_\m J(\m_k,\ta_k) \qquad \text{and} \qquad \nabla_\m \phi_k(\m_{k-1})=\nabla_\m J(\m_{k-1},\ta_k).
\end{align*}
 The former is satisfied by construction. The latter gives 
\[
\widehat{\vec{H}}_k(\m_k-\m_{k-1}) = \nabla_\m J(\m_k,\ta_k) - \nabla_\m J(\m_{k-1},\ta_k) .
\]
The gradient $\nabla_\m J(\m_k,\ta_k)$ is computed as a by-product of the 
pseudo-time continuation algorithm; see Appendix~\ref{appendix:computation}.  However,
computing $ \nabla_\m J(\m_{k-1},\ta_k)$ would require an additional forward and
adjoint solve. To avoid this, we consider a first-order Taylor expansion in
$\ta$,

\[
\nabla_\m J(\m_{k-1},\ta_k) \approx \nabla_\m J(\m_{k-1},\ta_{k-1}) + \Delta t \B(\m_{k-1},\ta_{k-1}) \Delta \ta .
\]
This implies the approximate secant condition
\begin{align}
\label{equ:secant_eqn_0}
\widehat{\vec{H}}_k \vec{z}_k = \vec{y}_k
\end{align}
where $\vec{z}_k = \m_k-\m_{k-1}$ and
\[
\vec{y}_k= \nabla_\m J(\m_k,\ta_k) - \nabla_\m J(\m_{k-1},\ta_{k-1}) - \Delta t \B(\m_{k-1},\ta_{k-1}) \Delta \ta.
\]
Note that $\vec{y}_k$ is defined 
in terms of vectors that are already computed.

Mimicking the BFGS update~\cite[Chapter 6]{NocedalWright06}, rather than requiring~\eqref{equ:secant_eqn_0}, we multiply~\eqref{equ:secant_eqn_0} by $\E_k=\widehat{\vec{H}}_k^{-1}$ and enforce the secant condition
\begin{align}
\label{equ:secant_eqn}
\E_k \vec{y}_k = \vec{z}_k .
\end{align}
We seek to determine a symmetric update matrix $\Delta \E_{k-1}$ such that $\E_k
= \E_{k-1} + \Delta \E_{k-1}$ and $\E_k$ satisfies~\eqref{equ:secant_eqn}. 
To determine an optimal $\Delta \E_{k-1}$, we consider 
\begin{align}
\label{equ:bfgs_opt_prob}
& \min_{\Delta \E_{k-1} \in \R^{n \times n}} \| \Delta \E_{k-1} \|^2 \\
& \Delta \E_{k-1} \vec{y}_k = \vec{x}_k, \nonumber \\
& \Delta \E_{k-1} = \Delta \E_{k-1}^\top, \nonumber
\end{align}
where $\vec{x}_k = \vec{z}_k - \E_{k-1} \vec{y}_k$ and $\| \cdot \|$ 
is suitably defined matrix norm. 
As detailed in~\cite[Chapter 6]{NocedalWright06}, \eqref{equ:bfgs_opt_prob}
admits an analytic solution, which leads to the update 
\begin{align}
\label{equ:bfgs_update}
\E_k = \left( \II - \rho_k \vec{z}_k \vec{y}_k^\top \right) \E_{k-1} \left( \II - \rho_k \vec{y}_k \vec{z}_k^\top \right) + \rho_k \vec{z}_k \vec{z}_k^\top,
\end{align}
where $\rho_k \coloneq 1/(\vec{y}_k^\top\vec{z}_k).$ Here, we require that the \emph{curvature condition} $\rho_k > 0$ 
holds. If it is not satisfied, the update~\eqref{equ:bfgs_update} is omitted. However, empirical observations
suggest that $\rho_k > 0$ is satisfied as long as the pseudo-time continuation procedure is faithfully tracking a minimizer 
about which $\nabla_{\m,\m} J$ is positive definite.
Note that the form of~\eqref{equ:bfgs_update} is identical to the classical BFGS update.
The difference is 
in the definition of $\vec{y}_k$.

\subsection{Block Quasi-Newton Update}

The parametric quasi-Newton update introduced in the previous subsection
incorporates the time evolution of $\ta(t)$ directly into the preconditioner. 
This alone is beneficial for preconditioning the CG iterations.  However, we can further
enhance the quasi-Newton update by using the Hessian-vector products computed
within PCG iterations. 
Specifically,
after executing
the PCG solve at time step $t_k$, we have a set of matrix-vector products $\vec{w}_i = \H(\m_k,\ta_k) \vec{p}_i$,
$i= 1,2,\ldots,\ell$ (see Line~11 of Algorithm~\ref{alg:pcg}), where we assume PCG has converged in $\ell$
iterations. 
The data 
$\{(\vec{w}_i, \vec{p}_i)\}_{i=1}^\ell$ can be used to enhance $\E_k$ 
before 
applying the parametric quasi-Newton
update~\eqref{equ:bfgs_update} to obtain $\E_{k+1}$. We discuss the details of this 
process next. 

Traditionally, quasi-Newton methods consider rank 1 or rank 2 
updates that satisfy a single secant equation. 
To incorporate the data 
\[
\vec{P}_k=\begin{pmatrix} \vec{p}_1 & \vec{p}_2 & \dots & \vec{p}_\ell \end{pmatrix}
\quad \text{and} \quad \vec{W}_k=\begin{pmatrix} \vec{w}_1 &
\vec{w}_2 & \dots & \vec{w}_\ell \end{pmatrix},
\]
we seek to enforce $\ell$ secant equations as stated in the following block form,
\begin{align} \label{eqn:block_secant}
	\widehat{\vec{H}}_k \vec{P}_k = \vec{W}_k. 
\end{align}
This amounts to requiring that the quasi-Newton Hessian approximation $\widehat{\vec{H}}_k$ agrees with the full Hessian 
on the Krylov space explored in the PCG iterations. 
As shown in~\cite{Schnabel_1983},
a symmetric matrix satisfying~\eqref{eqn:block_secant} exists
if and only if $\vec{P}_k^\top \vec{W}_k$ is symmetric. 
In the present setting, this condition is satisfied since
$\vec{P}_k^\top \vec{W}_k=\vec{P}_k^\top \H(\m_k,\ta_k) \vec{P}_k$.
We leverage this observation
and define a block BFGS update following the approach 
introduced in~\cite{gao_2018}. 
Specifically, 
given a preconditioner
$\E_k$ and PCG data $(\vec{P}_k, \vec{W}_k)$, the symmetric matrix
$\E_k^{\text{up}}$ that is closest to $\E_k$ and satisfies 
$\E_k^{\text{up}} \vec{W}_k = \vec{P}_k$ is
\begin{equation} \label{equ:block_bfgs}
\E_k^{\text{up}} = \left( \mathbf{I} - \vec{P}_k \vec{D}_k^{-1} \vec{W}_k^\top \right) 
						\E_k 
						\left( \vec{I} - \vec{W}_k \vec{D}_k^{-1} \vec{P}_k^\top \right) + 
						\vec{P}_k \vec{D}_k^{-1} \vec{P}_k^\top,
\end{equation}
where $\vec{D}_k=\vec{P}_k^\top \vec{W}_k \in \R^{\ell \times \ell}$. 

Note that the PCG search directions $\{\vec{p}_i\}_{i=1}^\ell$ satisfy the
orthogonality property $\vec{p}_i^\top \vec{w}_j = 0$ for $i \ne j$.  Therefore, in
exact arithmetic, $\vec{D}_k$ is a diagonal matrix. 
Since the block BFGS update~\eqref{equ:block_bfgs} depends on
$\vec{D}_k^{-1}$, care must be taken if a diagonal entry $\vec{p}_i^\top \vec{w}_i$
is small. Following the filtering idea presented in~\cite{gao_2018}, we discard
vectors for which $\vec{p}_i^\top \vec{w}_i$ is small. This is done by specifying a
tolerance $\tau \ge 0$ and discarding any vector for which $\vec{p}_i^\top
\vec{w}_i < \tau \| \vec{p}_i \|^2$. If $\tau = 0$, we retain all vectors
and as $\tau \to \infty$ we discard all vectors and hence do not perform the
update. We used $\tau=10^{-6}$ as a default value in our numerical results. 

However, in practice the aforementioned orthogonality of PCG search directions is
not preserved numerically. This loss of orthogonality, which is a well-known
phenomenon, causes the matrix $\vec{D}_k$ to be non-diagonal. Nonetheless, $\vec{D}_k$ 
is symmetric up to numerical precision. Therefore, in our implementation, we compute the 
eigenvalue decomposition $\vec{D}_k = \vec{V}_k \vec{\Lambda}_k \vec{V}_k^\top$,
and modify~\eqref{equ:block_bfgs} by replacing $\vec{P}_k$, $\vec{W}_k$, and $\vec{D}_k$ with $\vec{P}_k' = \vec{P}_k \vec{V}_k$, $\vec{W}_k'=\vec{W}_k\vec{V}_k$, and $\vec{D}_k'=\vec{P}_k'^\top \vec{W}_k'$, respectively. 
Since the matrix $\vec{D}_k$ is typically low-dimensional, the eigenvalue decomposition may be computed
with standard dense linear algebra routines.
Rotating $\vec{P}_k$ and $\vec{W}_k$ with the eigenvectors $\vec{V}_k$ is equivalent to enforcing the secant condition $\E_k^{\text{up}} \vec{W}_k' = \vec{P}_k'$ and leads to the numerically diagonal matrix
$\vec{D}_k'=\vec{P}_k'^\top \vec{W}_k'=\vec{V}_k^\top \vec{P}_k^\top \vec{W}_k \vec{V}_k=\vec{V}_k^\top \vec{V}_k \vec{\Lambda}_k \vec{V}_k^\top \vec{V}_k=\vec{\Lambda}_k$.

\subsection{Algorithm Summary}\label{sec:algorithm_summary}
Algorithm~\ref{alg:precond_pt_contin} summarizes preconditioned forward Euler
pseudo-time continuation and highlights the preconditioner adaptation. In
the loop over time steps $k=0,1,2\dots,N-1$, we first compute a Forward
Euler time step to produce $\m_{k+1}^\text{pred} \approx \ms(\ta_{k+1})$. Next, we use
a block quasi-Newton update to improve the time step $t_k$ preconditioner using
$\H(\m_k,\ta_k)$ matrix-vector products, and subsequently produce a time step
$t_{k+1}$ preconditioner via a parametric quasi-Newton update. The corrector
step executes a Newton optimization iteration at time step $t_{k+1}$, and
subsequently updates the time step $t_{k+1}$ preconditioner via a block
quasi-Newton update using $\H(\m_{k+1}^\text{pred},\ta_{k+1})$ matrix-vector products. 
The tolerance satisfaction step may not be required as the predictor/corrector
scheme typically predicts $\ms(\ta_{k+1})$ effectively. However, if the
time step size $\Delta t$ is large, the predictor/corrector scheme may be
inaccurate. The subsequent tolerance satisfaction step executes Newton
optimization iterations until convergence. The PCG tolerance $\epsilon_\text{CG}$ 
may be reduced in the tolerance satisfaction step to ensure Newton convergence properties. 
However, in our numerical experiments, 
tightening this tolerance yielded no benefits.

The initial preconditioner rank $r_\text{init} \ge 0$ specifies $\E_0$ as: (i) the
inverse regularization Hessian if $r_\text{init}=0$, or (ii) an approximate inverse
based on a rank $r_\text{init}$ approximation of $\vec{R}^{\frac{1}{2}}
\HM(\m_0,\ta_0)\vec{R}^{\frac{1}{2}}$ if $r_\text{init}>0$.  Thanks to the block
quasi-Newton update, our proposed approach has a self-correcting property. If
the preconditioner is a poor approximation of the inverse Hessian, then more PCG
iterations will be required. These iterations provide data to enhance the
preconditioner so that it better approximates the inverse Hessian and hence
enhances the future preconditioners. The maximum update rank $r_\text{update}
\ge 0$, along with the tolerance $\tau$, controls the rank of the block
Quasi-Newton update. In practice, $r_\text{update}$ is
specified based on storage limitations.

We can also consider a version of Algorithm~\ref{alg:precond_pt_contin} that uses modified Euler in place of forward Euler in the predictor step. This involves two
Hessian inversions and correspondingly two preconditioner updates within the
predictor step. However, this additional cost may be offset by the higher-order
predictor method requiring fewer time steps, PCG iterations, and/or tolerance
satisfaction iterations.

As a final remark, note that we check the norm of the gradient before
executing the corrector step.  If the optimality tolerance is satisfied, we
forego the corrector step.  This can lead to considerable computational savings.
We revisit this issue in our numerical experiments in 
Section~\ref{sec:numerics_icesheet}.

\begin{algorithm}
\caption{Preconditioned forward Euler Pseudo-time Continuation}
\begin{algorithmic}
\STATE \textbf{Input: } $N>0$, $r_\text{init} \ge 0$, $r_\text{update} \ge 0$, $\epsilon_{\text{CG}}>0$
\STATE Initialize $\E_0$, a rank $r_\text{init}$ preconditioner for $\H(\m_0,\ta_0)$
\FOR{$k=0,1,\dots,N-1$}
\STATE \hspace{-1em} \noindent\rule{\linewidth}{0.2mm}
\STATE \underline{\textbf{Predictor step:}}
\STATE Compute $\vec{g}_k^\text{pred} = \Delta t \B(\m_k,\ta_k) \Delta \ta$
\STATE Compute $[\vec{s}_k^\text{pred},\vec{P}_k,\vec{W}_k] = \text{PCG}(\vec{g}_k^\text{pred},\H(\m_k,\ta_k),\E_k,\epsilon_{\text{CG}})$ 
\hfill\COMMENT{Alg.~\ref{alg:pcg}}\hspace{-1em}
\STATE Set $\m_{k+1}^\text{pred} = \m_k -  \vec{s}_k^\text{pred}$
\STATE Set $\E_k = \text{Block\_QN\_Update}(\E_k,\vec{P}_k,\W_k,r_\text{update})$
\hfill\COMMENT{see~\eqref{equ:block_bfgs}}\hspace{-1em}
\STATE \hspace{-1em} \noindent\rule{\linewidth}{0.2mm}
\STATE \underline{\textbf{Corrector step:}}
\STATE Compute $\vec{g}_{k+1} = \nabla_\m J(\m_{k+1}^\text{pred},\ta_{k+1})$
\STATE Set $\E_{k+1} = \text{Parametric\_QN\_Update}(\E_k,\vec{g}_k^\text{pred},\vec{g}_k,\vec{g}_{k+1},\m_k,\m_{k+1}^\text{pred})$
\hfill\COMMENT{see~\eqref{equ:bfgs_update}}\hspace{-1em}
\STATE Compute $[\vec{s}_{k+1},\vec{P}_{k+1},\W_{k+1}] \hspace{-.75mm}=\hspace{-.75mm} \text{PCG}(\vec{g}_{k+1},\H(\m_{k+1},\ta_{k+1}),\E_{k+1},\epsilon_{\text{CG}})$  
\hspace{-1.8mm} \hfill\COMMENT{Alg.~\ref{alg:pcg}} \hspace{-1.3em}
\STATE Set $\m_{k+1} = \m_{k+1}^\text{pred} - \vec{s}_{k+1}$
\STATE Set $\E_{k+1} = \text{Block\_QN\_Update}(\E_{k+1},\vec{P}_{k+1},\W_{k+1},r_\text{update})$
\hfill\COMMENT{see~\eqref{equ:block_bfgs}}  \hspace{-1.3em}
\STATE \hspace{-1em} \noindent\rule{\linewidth}{0.2mm}
\STATE \underline{\textbf{Tolerance satisfaction step:}}
\STATE Compute $\vec{g}_{k+1} = \nabla_\m J(\m_{k+1},\ta_{k+1})$
\WHILE{$\| \vec{g}_{k+1} \| \ge$ tolerance}
\STATE Compute $\vec{s}_{k+1} = \text{PCG}(\vec{g}_{k+1},\H(\m_{k+1},\ta_{k+1}),\E_{k+1},\epsilon_{\text{CG}})$
\hfill\COMMENT{Alg.~\ref{alg:pcg}} \hspace{-1em}
\STATE Set $\m_{k+1} = \m_{k+1} - \vec{s}_{k+1}$
\STATE Compute $\vec{g}_{k+1} = \nabla_\m J(\m_{k+1},\ta_{k+1})$
\ENDWHILE
\STATE \hspace{-1em} \noindent\rule{\linewidth}{0.2mm}
\ENDFOR
\end{algorithmic} 
\label{alg:precond_pt_contin} 
\end{algorithm}

\section{Numerical results} \label{sec:numerical_results}
In this section, we present numerical results for two nonlinear inverse
problems. Section~\ref{ssec:diff_example} presents a study on the effect of
algorithm parameters in our preconditioned pseudo-time continuation algorithm.
These results focus on providing intuition about the algorithm.
Section~\ref{sec:numerics_icesheet} demonstrates the use of preconditioned
pseudo-time continuation for a challenging ice sheet bedrock inversion problem.
The results in this section highlight the utility of our proposed approach for 
nonlinear and computationally intensive inverse problems. 

\subsection{Diffusion coefficient estimation} \label{ssec:diff_example}
We consider the Poisson equation
\begin{alignat}{2} \label{equ:diff}
 -\nabla \left( e^m \nabla u \right) &= h_{\ta} \qquad  &\quad \text{ on } \Omega,   \\
 e^m \nabla u \cdot \vec{n} &= 0 &\quad \text{ on } \Gamma_1 \cup \Gamma_3, \nonumber \\
 u &= g_2 &\quad \text{ on } \Gamma_2, \nonumber \\
 u &= g_4 &\quad \text{ on } \Gamma_4, \nonumber 
\end{alignat}
defined on the domain $\Omega = (0,1)^2$ whose boundary consists of $\Gamma_1 =
\{0 \} \times (0,1)$, $\Gamma_2 = (0,1) \times \{0\} $, $\Gamma_3 = \{1\} \times
(0,1)$, and $\Gamma_4 = (0,1) \times \{1\}$. Here, $\vec{n}$ denotes the outward-facing
 normal vector to the boundary. In this problem $h_{\ta}$ is a prescribed (parameterized) source
term and $g_2$ and $g_4$ are predefined functions. We use $g_2(x) = \cos(4 \pi x)^2$
and $g_4(y) = \sin(2 \pi x)^2$. Consider  
the inverse problem of estimating the spatially distributed 
log-diffusion coefficient $m:\Omega \to \R$.

In this problem, the source term $h_{\ta}$ is parameterized by 
an auxiliary parameter vector 
$\ta \in \R^9$ as follows:
\begin{align*}
h_{\ta}(x, y) = 100.0 \sum\limits_{i=1}^3 \sum\limits_{j=1}^3 \theta_{i,j} \frac{\sin(2 \pi i x) \sin(2 \pi j y)}{i j} .
\end{align*}

\begin{figure}[h]
\centering
    \includegraphics[width=0.40\textwidth]{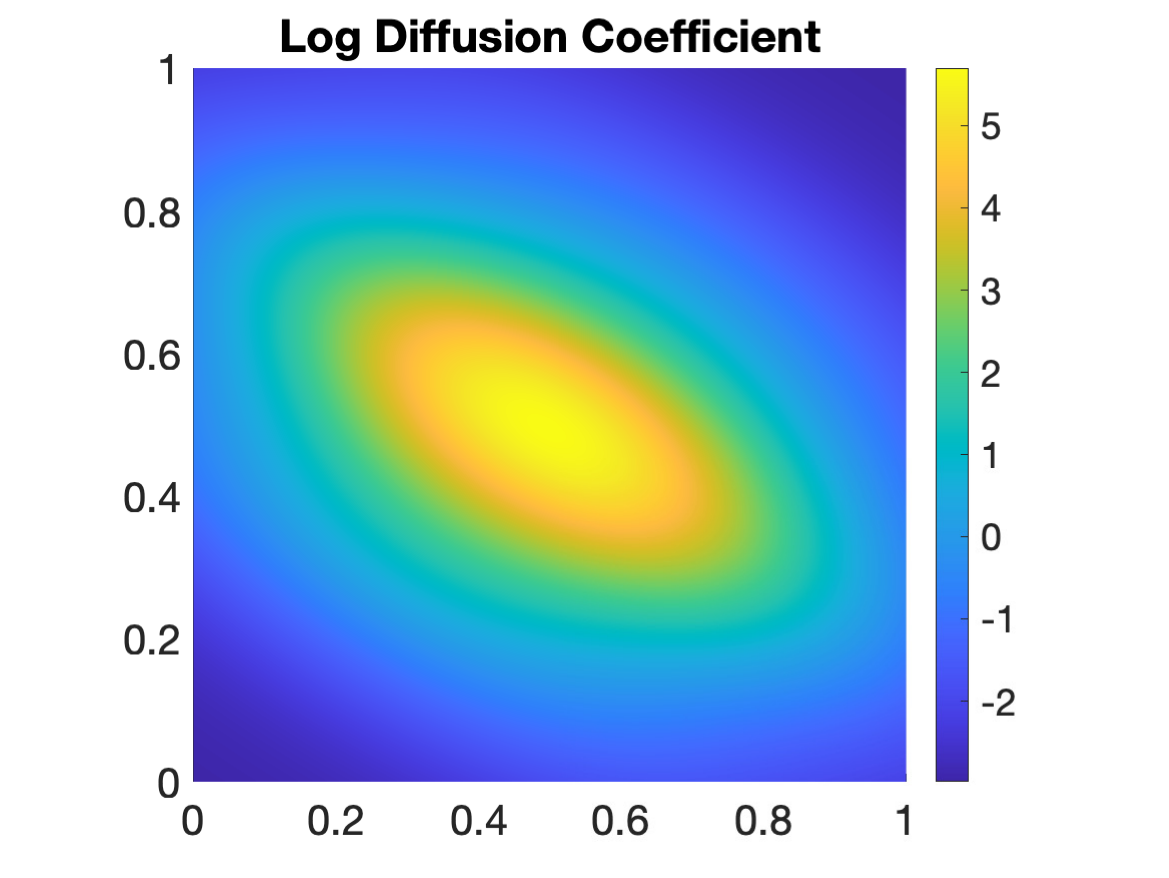}
  \includegraphics[width=0.40\textwidth]{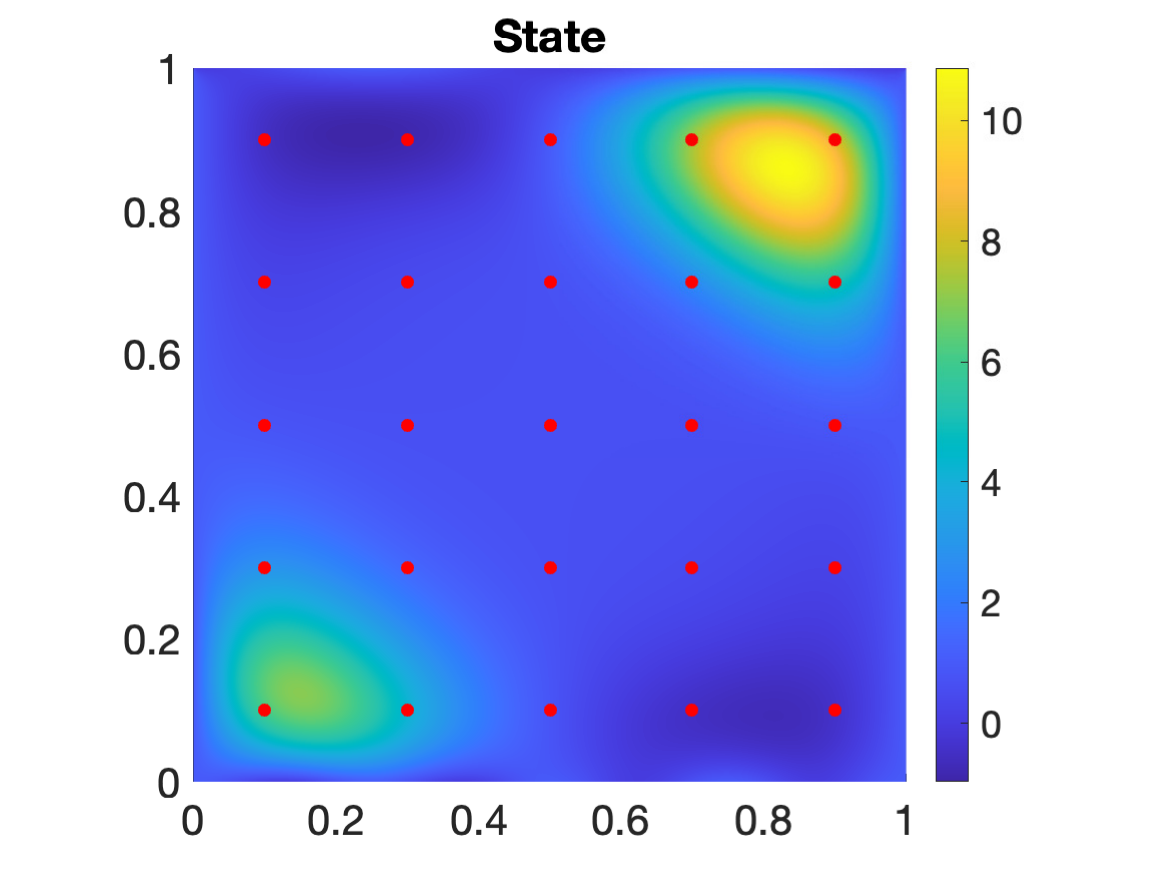}
    \caption{Left: log diffusion coefficient used to generate synthetic data; right: the state solution with observation locations denoted by red dots.}
  \label{fig:state_and_diff_coeff}
\end{figure}

To solve the inverse problem, we discretize~\eqref{equ:diff} on a $50 \times 50$
rectangular mesh. Hence, the discretized state variable and inversion parameter have
dimension $51^2=2601$.  In the present study, we use the ground-truth parameter
$m^\dagger$ depicted in the left panel of Figure~\ref{fig:state_and_diff_coeff}
and use a ground-truth auxiliary parameter of $\ta^\dagger$, which will be 
specified below. Subsequently, we generate synthetic data by 
solving~\eqref{equ:diff} with $(m, \ta) = (m^\dagger, \ta^\dagger)$ and
extracting solution values on a grid of measurement points depicted in the right
panel of Figure~\ref{fig:state_and_diff_coeff}. 
To avoid the inverse crime, the data is generated by solving the forward problem
on a $200 \times 200$ rectangular mesh and contaminated with $2\%$ Gaussian noise.  

We pose an inverse problem in the form of~\eqref{equ:optim} where $\mathbf{F}$
corresponds to the composition of the (discretized) solution operator of~\eqref{equ:diff} and
the observation operator maps the state solution to sparse observations.
Furthermore, $\mathbf{R}$ corresponds to the square of the discretized elliptic
operator $\gamma(-\kappa \Delta + \mathcal I)$, with constants $\kappa =10$ and $\gamma =10^{-3}$, where $\Delta$ and $\mathcal I$
correspond to the Laplacian and identity operators, respectively.  This choice
is common in practice as it corresponds to using the squared inverse elliptic
operator to define the prior covariance in the Bayesian formulation of the
inverse problem; see~\cite{Bui-ThanhGhattasMartinEtAl13}.

When solving the inverse problem, we use a nominal source function
parameterization $\overline{\ta}$ rather than the ground-truth auxiliary
parameter $\ta^\dagger$ that was used to generate the data.  Subsequently, 
we use pseudo-time continuation to compute the solution of the inverse problem when
$\ta = \tat$, where $\tat$ is an auxiliary parameter vector that approximates
$\ta^\dagger$.

\subsubsection*{Varying the number of time steps $N$}

\begin{figure}[h]
\centering
  \includegraphics[width=0.32\textwidth]{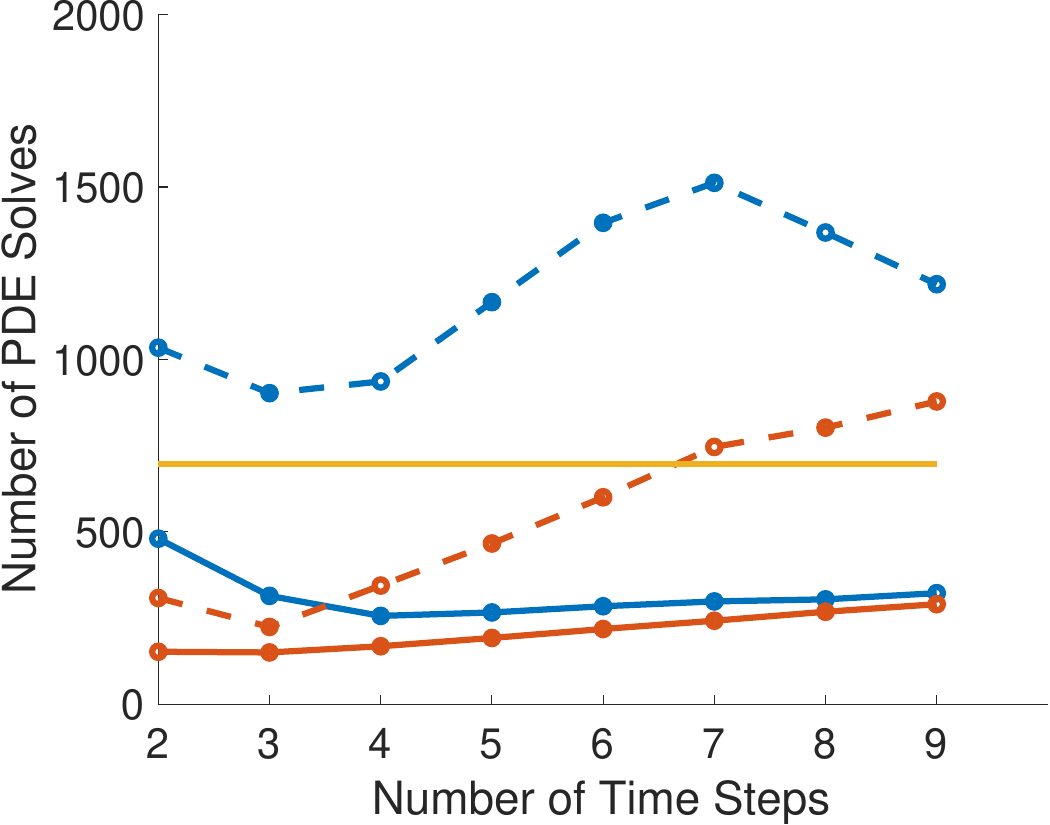}
    \includegraphics[width=0.32\textwidth]{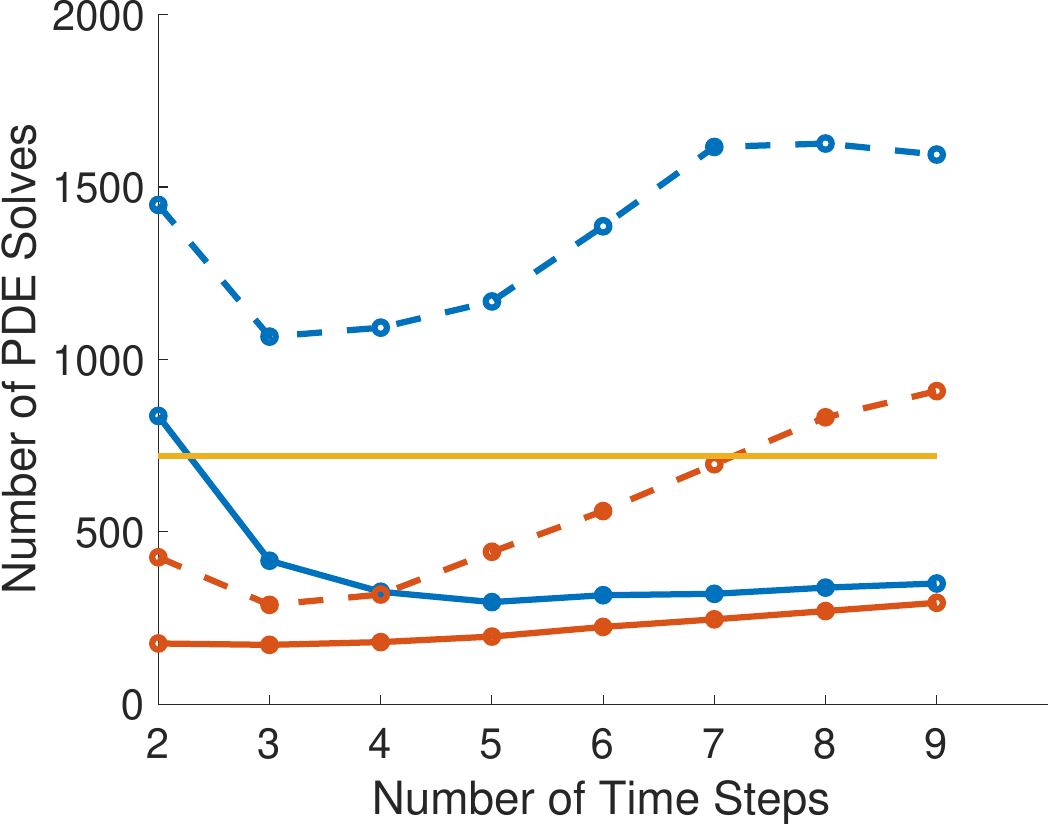}
      \includegraphics[width=0.32\textwidth]{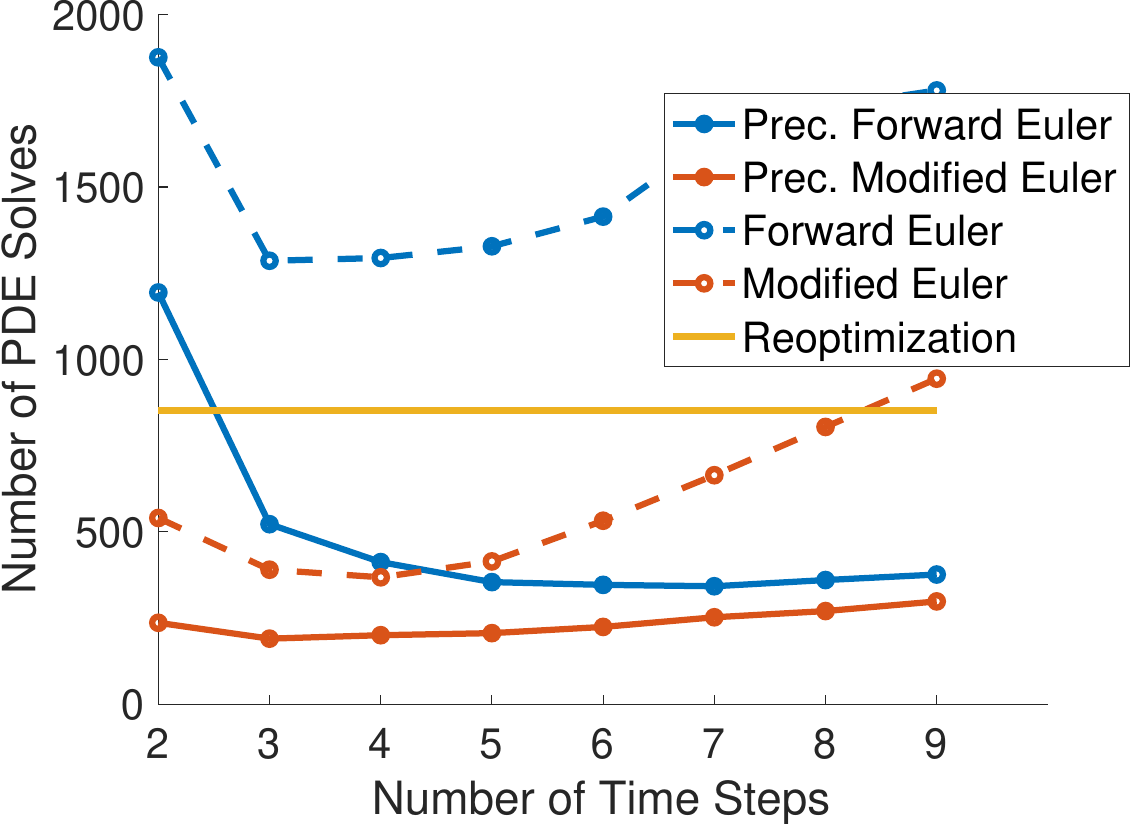}
    \caption{Comparison of Forward and modified Euler time-stepping for three different $\ta$ perturbation magnitudes $\alpha=1.0$ (left), $\alpha=1.25$ (center), and $\alpha=1.5$ (right).}
  \label{fig:cost_comparison_vary_timesteps}
\end{figure}

In our first numerical experiment, we compare the use of forward Euler and
modified Euler time-stepping algorithms for the predictor. We fix $\overline{\ta}=\e_1$, the
canonical basis vector whose first entry is $1$ and all others are $0$, and
consider three perturbations of varying magnitudes. Specifically, we consider
$\tat = \overline{\ta} + \alpha \e$, where $\e \in \R^9$ is the vector
whose entries are all $1$, and $\alpha \in \{1.0, 1.25, 1.5\}$.
Note that $\ta^\dagger=\overline{\ta}+1.3 \e$. For each scenario, we
compare the performance of five potential approaches. Specifically, we
establish a performance baseline by solving the perturbed optimization
problem~\eqref{equ:optim} with $\ta=\tat$ and the initial condition given
by the optimal solution when $\ta=\overline{\ta}$. We refer to this as ``re-optimization".
 In addition, we execute the pseudo-time continuation procedure
in four different configurations: using either forward Euler or modified Euler for the predictor steps, 
and using either the inverse regularization Hessian as the preconditioner for all solves or using
our proposed adaptive preconditioner. In the four pseudo-time continuation configurations, we
fix the PCG tolerance to $\epsilon_\text{CG} = 10^{-4}$ and vary
the numbers of time-steps; specifically $N \in \{2,3,\dots,9\}$. 
For our proposed adaptive preconditioner, we use the inverse regularization Hessian as the initial preconditioner and set the
maximum rank of each block Quasi-Newton update to be $r_\text{update} = 20$. 
The results are shown in Figure~\ref{fig:cost_comparison_vary_timesteps}, where we
measure the computational cost via the total number of PDE solves required to
achieve a gradient norm tolerance of $10^{-8}$.

We make several observations from Figure~\ref{fig:cost_comparison_vary_timesteps}.
\begin{itemize}
\item A modified Euler predictor is generally superior to a forward Euler 
predictor; however, the difference in performance is less pronounced
when the $\ta$ perturbation magnitude is small and the number
of time steps is large.

\item Using a fixed preconditioner results in considerable variability in computational cost as the number of time steps varies. In contrast,
using the adaptive preconditioner results in small changes in computational cost as the number of time steps varies.

\item Taking too few time steps, particularly with a forward Euler predictor step, may result in increased computational cost
  as a result of requiring multiple iterations to satisfy the gradient tolerance.
  
\item The computational cost of the pseudo-time continuation algorithms grows as the 
magnitude of the $\ta$ perturbation increases.

\item The computational cost of re-optimization is always larger than pseudo-time continuation
with a modified Euler predictor and an adaptive preconditioner. 

\item The computational cost of re-optimization is always less than pseudo-time continuation
with a forward Euler predictor and a fixed preconditioner. This is because re-optimization
takes less costly steps (fewer Hessian-vector products per iteration)
when it is initially far from the minimizer.

\end{itemize}

These observations highlight the benefit of using a modified Euler predictor with the adaptive preconditioner. 
However, the results above correspond to a fixed PCG
tolerance and preconditioner update rank. In what follows, we vary these algorithm parameters while
fixing the number of time steps ($N=3$) and the perturbation magnitude ($\alpha=1.25$).

\subsubsection*{Varying the PCG tolerance $\epsilon_\text{CG}$}

\begin{figure}[h]
\centering
  \includegraphics[width=0.40\textwidth]{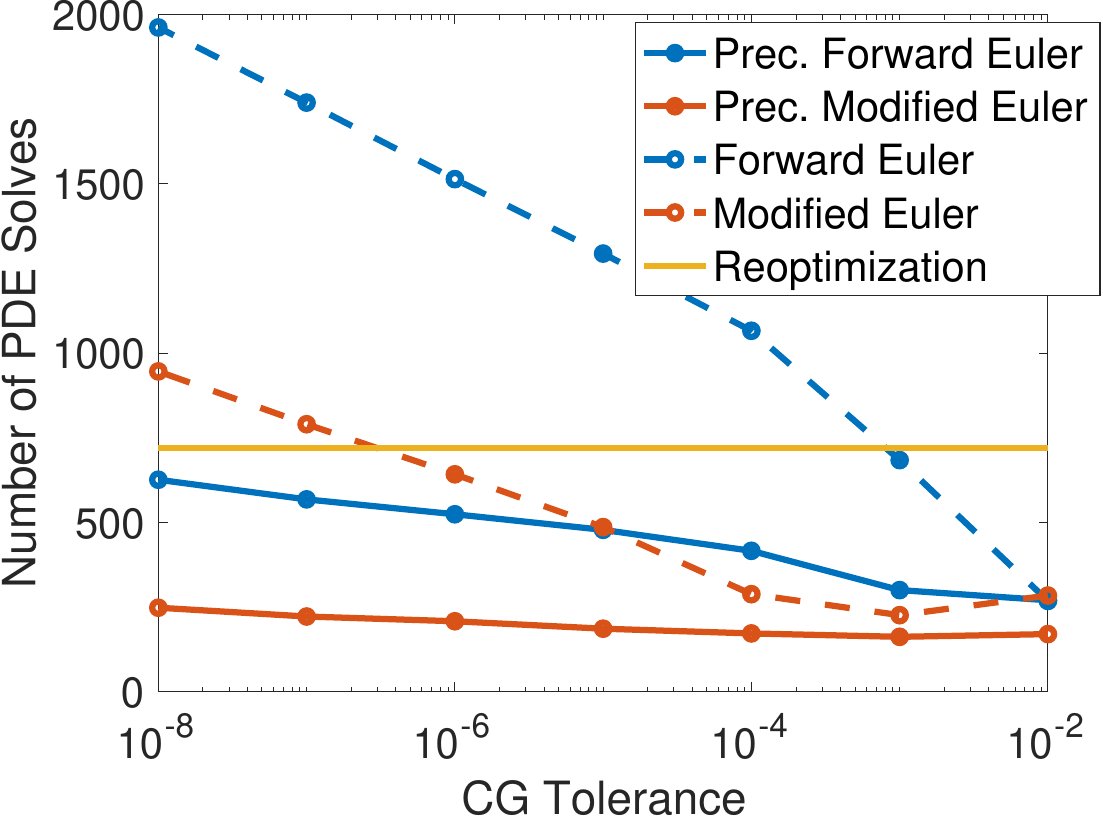}
    \caption{Comparison of Forward and modified Euler for various PCG tolerances.}
  \label{fig:cost_comparison_vary_tolerance}
\end{figure}

In our second numerical experiment, we consider the effect of changing the PCG
tolerance. We set the initial
preconditioner as the inverse regularization Hessian and the maximum rank
of the block Quasi-Newton update to be $r_\text{update} = 20$.
Figure~\ref{fig:cost_comparison_vary_tolerance} displays the computational cost
as a function of the PCG tolerance.  As before, the cost is measured in terms of the total
number of PDE solves required to achieve a gradient norm tolerance of $10^{-8}$.

We observe that the optimal PCG tolerance in this example is around 
$\epsilon_\text{CG} = 10^{-2}$. In this case, modified Euler with the adaptive preconditioner
is slightly superior to the other three pseudo-time continuation approaches, which all have
comparable computational costs. We observed that the tolerance $\epsilon_\text{CG} = 10^{-1}$ 
could be effective, but the time stepping procedure may not converge if the loose tolerance
leads to ineffective predictor steps that fail to track the Newton basin of attraction. 
Because of this lack of robustness, we omitted 
$\epsilon_\text{CG} = 10^{-1}$ from the plot. The success of our pseudo-time continuation with
such a loose PCG tolerance is explained through a combination of several factors.
First, it is unsurprising that a loose tolerance is beneficial in the predictor
step as the goal is to account for variations in $\ta$ and get sufficiently
close to the minimizer that the corrector step will be successful.  However,
surprisingly, a loose tolerance on the corrector step is superior in this
example.  Although additional tolerance satisfaction steps were required, the
total cost was lowest when using a loose tolerance. 

It is noteworthy that using
a fixed preconditioner results in a significant increase in computational cost as the the tolerance
is tightened. This is because the number of Hessian-vector products required in the PCG solve increases.
The adaptive preconditioner helps to mitigate this increased cost. However, in our later 
numerical example, we will observe that the computational cost using the adaptive preconditioner
can vary nontrivially with changes in the PCG tolerance.

\subsubsection*{Varying the preconditioner rank parameters $r_\text{init}$ 
and $r_\text{update}$} 
In our final numerical experiment, we focus on the performance of pseudo-time continuation
with the modified Euler predictor, as $r_\text{init}$ and $r_\text{update}$ are varied. In
particular, we initialize the preconditioner using a rank $r_\text{init}$
approximation of the regularization-preconditioned data-misfit Hessian and set
the maximum block Quasi-Newton update to be of rank $r_\text{update}$. We
consider all pairs $(r_\text{init},r_\text{update}) \in \{0,5,10,15,20\}^2$. 
For each $(r_\text{init},r_\text{update})$, we execute the pseudo-time
continuation procedure with a modified Euler predictor, $N=3$ time steps, and a
PCG tolerance of $\epsilon_\text{CG} = 10^{-2}$. 
Figure~\ref{fig:cost_comparison_vary_tolerance} displays the computational cost
(left panel) and the number of vectors stored (right panel) for various
rank parameters.

\begin{figure}[h]
\centering
  \includegraphics[width=0.40\textwidth]{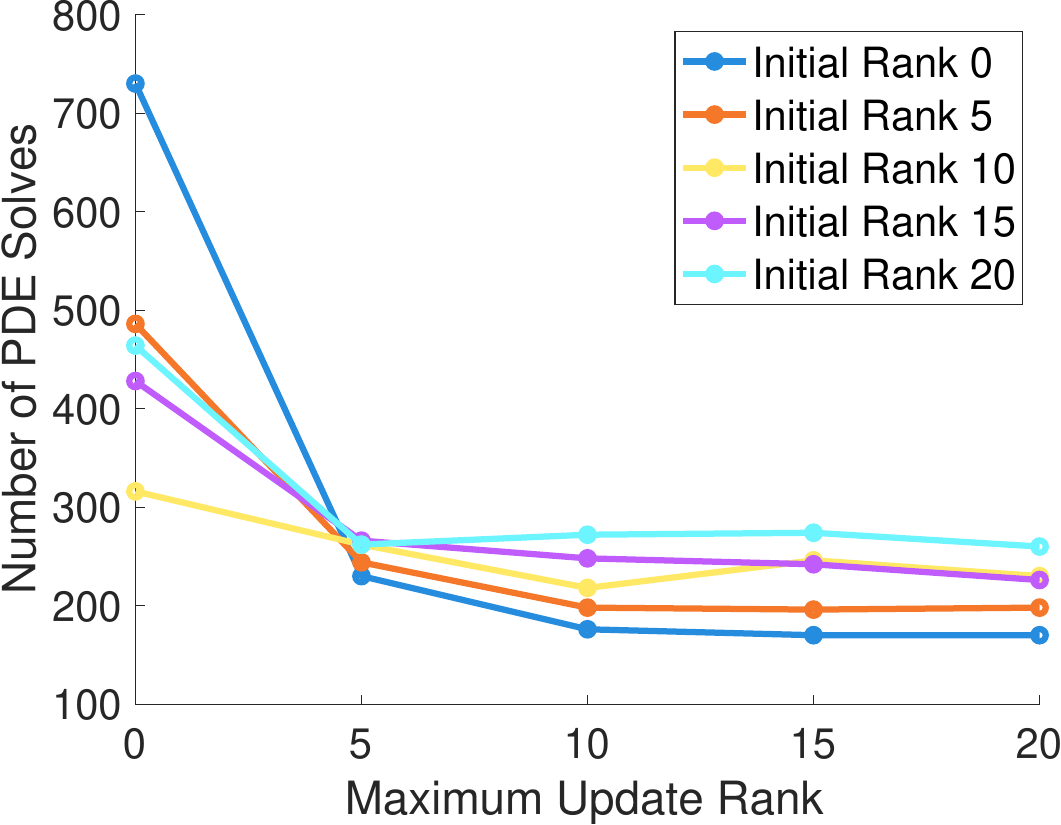}
    \includegraphics[width=0.40\textwidth]{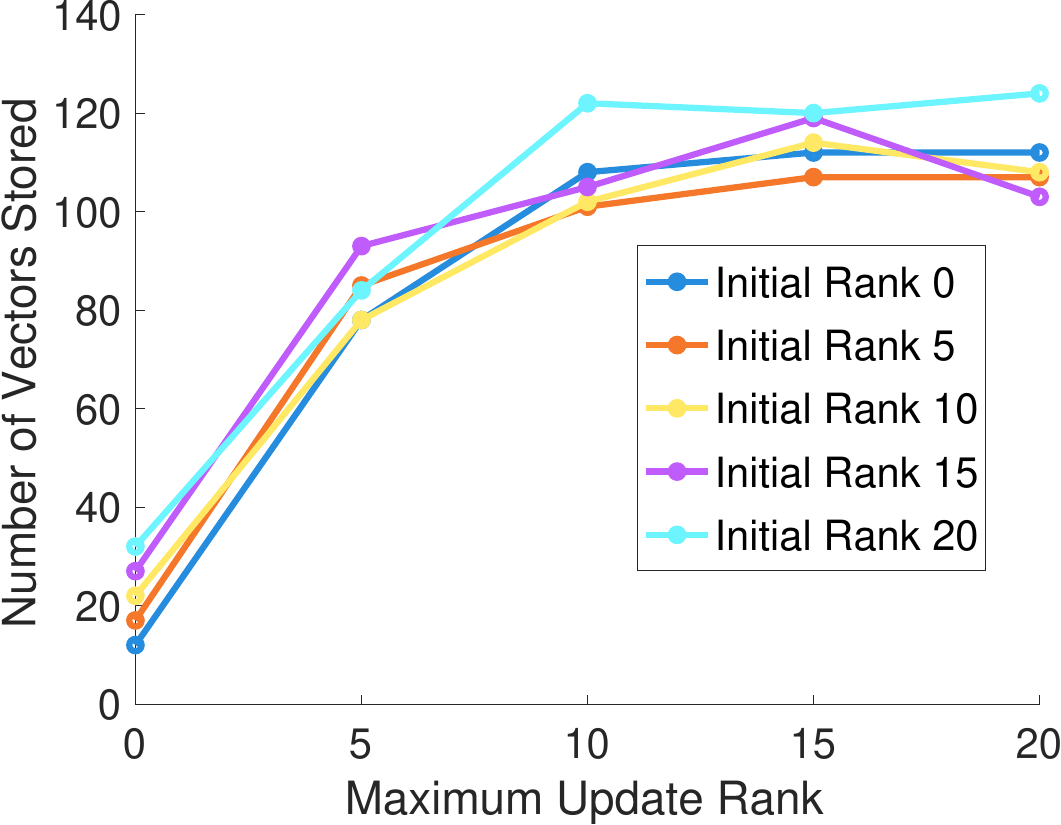}
    \caption{The cost of the pseudo-time continuation procedure for various
    preconditioner rank parameters. The computational cost in number of PDE solves is
    shown in the left panel and the number of vectors stored is shown in the
    right panel. Each curve corresponds to a rank $r_\text{init}$ used to
    initialize the preconditioner via a low-rank approximation of the
    regularization-preconditioned data-misfit Hessian.}

  \label{fig:cost_comparison_vary_rank}
\end{figure}

We make several observations from Figure~\ref{fig:cost_comparison_vary_rank}.
\begin{itemize}
\item There is a general trend of decreasing cost with increasing the
block Quasi-Newton update rank $r_\text{update}$, although diminishing returns
are observed after $r_\text{update}=5$.

\item There is a general trend of increasing storage with increasing the block
Quasi-Newton update rank $r_\text{update}$; however, the rate of storage growth
slows for larger $r_\text{update}$.

\item Setting $r_\text{update}=0$ (i.e., not using the block Quasi-Newton update)
results in a considerable increase in computational cost.

\item Investing in the initial preconditioner (i.e., taking a larger $r_\text{init}$) is advantageous when $r_\text{update}=0$, but has little benefit for larger $r_\text{update}$.
\end{itemize}

\subsubsection*{Scalability and mesh-independence}
The proposed computational framework exploits specific problem structures that
are independent of the discretization. Namely, the spectral properties of the
Hessian operator depend primarily on the regularity of the
parameter-to-observable map as well as smoothing properties of the
regularization operator. The Newton--Krylov methods exploit these structures to 
achieve a nearly mesh-independent behavior, as observed in various studies; see,
e.g.,~\cite{GhattasWilcox21} for discussions and references. 

To assess the mesh-independence properties of our proposed framework, we fixed
all algorithm parameters and varied the mesh resolution.  We observed negligible
variation in the computational cost as the mesh resolution varied, thus
indicating that the algorithm's convergence behavior depends on the inherent
spectral properties of the second-order derivative operators rather than the
dimension of the discretized problem. The results are omitted for brevity. 

\subsection{Ice sheet bedrock inversion}\label{sec:numerics_icesheet}

We consider estimating the bedrock topography on the bottom of the
Greenland Ice Sheet. For simplicity, we focus on a $550 \times 450$ km
region of the ice sheet (denoted as $\Omega$) and use the shallow ice
approximation \cite{morland_1980,hutter_1983} with the isothermal assumption \cite{bueler_2005}. In what follows, we summarize the inverse problem and use it as a testbed to analyze computational aspects of our pseudo-time continuation approach. We refer the reader to~\cite{Reese_2024}, where the inverse problem and its sensitivities were studied, for a more detailed discussion. 

The forward model is
\begin{alignat}{2} \label{eqn:sia_pde}
 \frac{\partial s}{\partial t} - \nabla \cdot (Q(s)\nabla s) =& h_{\ta} \qquad & \text{on } \Omega \times (0,T]  \\
\nabla s \cdot n = & 0 & \text{on } \partial \Omega \times (0,T] \nonumber \\
 s =&  s_0 & \text{on } \Omega \times \{0\} \nonumber \\
Q(s) =& e^{-\beta_{\ta}} \rho g (s-m)^2 + \frac{2A\rho^3 g^3}{5}(s-m)^5 \vert \vert \nabla s \vert \vert^2, \nonumber
\end{alignat}
where $s$ is the surface height of the ice sheet (relative to sea level), $Q(s)$ is a velocity field, $\beta_{\ta}$ is the log basal friction on the interface of the ice sheet and the bedrock beneath it, and $h_{\ta}$ is a forcing term modeling the accumulation/ablation of ice. Scalar parameters include: the density of ice $\rho = 910$ $(kg/m^3)$, acceleration due to gravity $g=9.81$ $(m/s^2)$, and a flow rate factor $A=10^{-16}$ $(Pa^{-3}/s)$. The forward model~\eqref{eqn:sia_pde} is transient. We integrate over the time interval from $0$ to $T=10$ years.

The inverse problem is to estimate the bedrock topography $m:\Omega \to \R$ using observations of the surface ice velocity, which is defined as a function of the PDE solution $s(m,\beta_{\ta},h_{\ta})$ via
\begin{equation}
    \vec{v}(m,\beta_{\ta},h_{\ta})=-\frac{1}{2} A \rho^3 g^3 (s(m,\beta_{\ta},h_{\ta})-m)^4 \vert \vert \nabla s(m,\beta_{\ta},h_{\ta}) \vert \vert^2 \nabla s(m,\beta_{\ta},h_{\ta}) \label{eqn:vel}.
\end{equation}

The log basal friction $\beta_{\ta}$ and forcing $h_{\ta}$ are spatially heterogeneous uncertain auxiliary parameters that we model as
\begin{equation*}
\beta_{\ta} = \overline{\beta}\delta_{1}(\ta) \qquad \text{and} \qquad h_{\ta} = \overline{h}\delta_{2}(\ta),
\end{equation*}
where 
\begin{equation*}
\delta_i(\ta) = 1 + 0.2\sum_{j=1}^{961} \theta_{((i-1)n+j)}\phi_j, \qquad i = 1,2,
\end{equation*}
are parameterized perturbations; $\{\phi_j\}_{j=1}^n$ are linear finite element
basis functions defined on a $31 \times 31$ rectangular mesh of the domain
$\Omega$. This implies that the auxiliary parameter dimension is
$(2)(31^2)=1922$. The nominal auxiliary parameters correspond to taking
$\overline{\ta} = \vec{0} \in \R^{1922}$ so that
$\beta_{\vec{0}}=\overline{\beta}$ and $h_{\vec{0}}=\overline{h}$. 

After discretization, the inverse problem takes the form of~\eqref{equ:optim} where $\vec{F}$, the parameter-to-observable map, is defined via composing the discretized solution operator for~\eqref{eqn:sia_pde} and velocity model~\eqref{eqn:vel}. In practice, observational data comes from satellite measurements whose spatial resolution is greater than the computational mesh resolution. Hence, observations are assumed to be available at all spatial locations. The regularization is defined as the negative log  of the prior covariance in a Bayesian formulation of the inverse problem, so that $\mathbf{R}$ corresponds to the square of the discretized elliptic
operator $\gamma(-\kappa \Delta + \mathcal I)$, with constants $\kappa =10^{-2}$ and $\gamma =9\times10^{-7}$.

Synthetic data is generated from the model using the ``true'' bedrock
topography, which is subsequently considered unknown as we solve the inverse
problem to reconstruct it. The bedrock, ice thickness, accumulation/ablation,
and log basal friction are taken from \cite{albany_felix}. To avoid an ``inverse
crime," we generate data by solving \eqref{eqn:sia_pde} on a $101\times 101$
mesh with $121$ time steps and interpolate the data onto a $71\times  71$ mesh
with $61$ time steps to solve the inverse problem. Additionally, $5\%$ Gaussian
noise is added to emulate the noise in satellite observations. Consequently, the
discretized state dimension, which corresponds to the observational data
dimension, is $(61)(71^2)=307,501$. The inversion parameter has dimension
$71^2=5041$. Note, however, that this is still an ill-posed inverse problem 
due to smoothing properties of the forward operator.  Furthermore, this is a
challenging inverse problem due to the complexity of the nonlinear forward model
and the auxiliary parameter uncertainty. 

\begin{figure}[h]
\centering
  \includegraphics[width=0.32\textwidth]{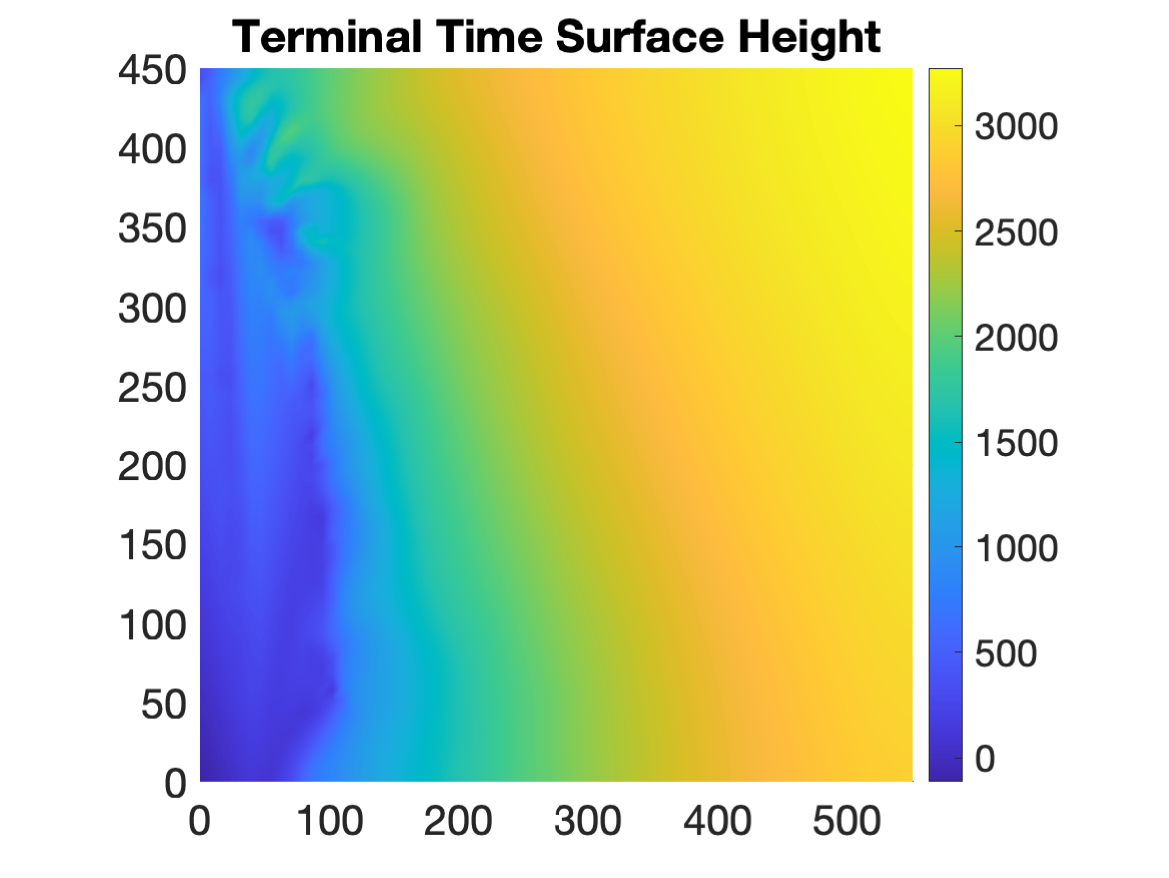}
    \includegraphics[width=0.32\textwidth]{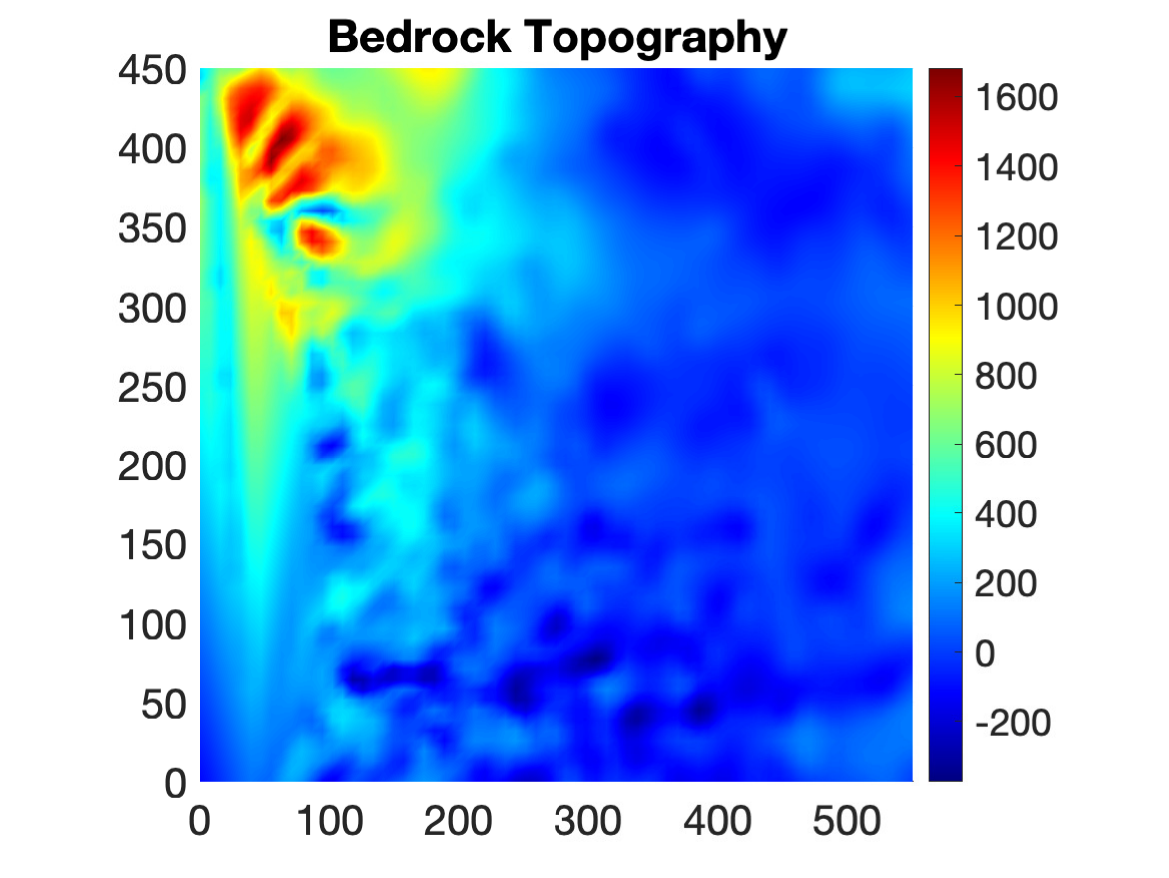}
       \includegraphics[width=0.32\textwidth]{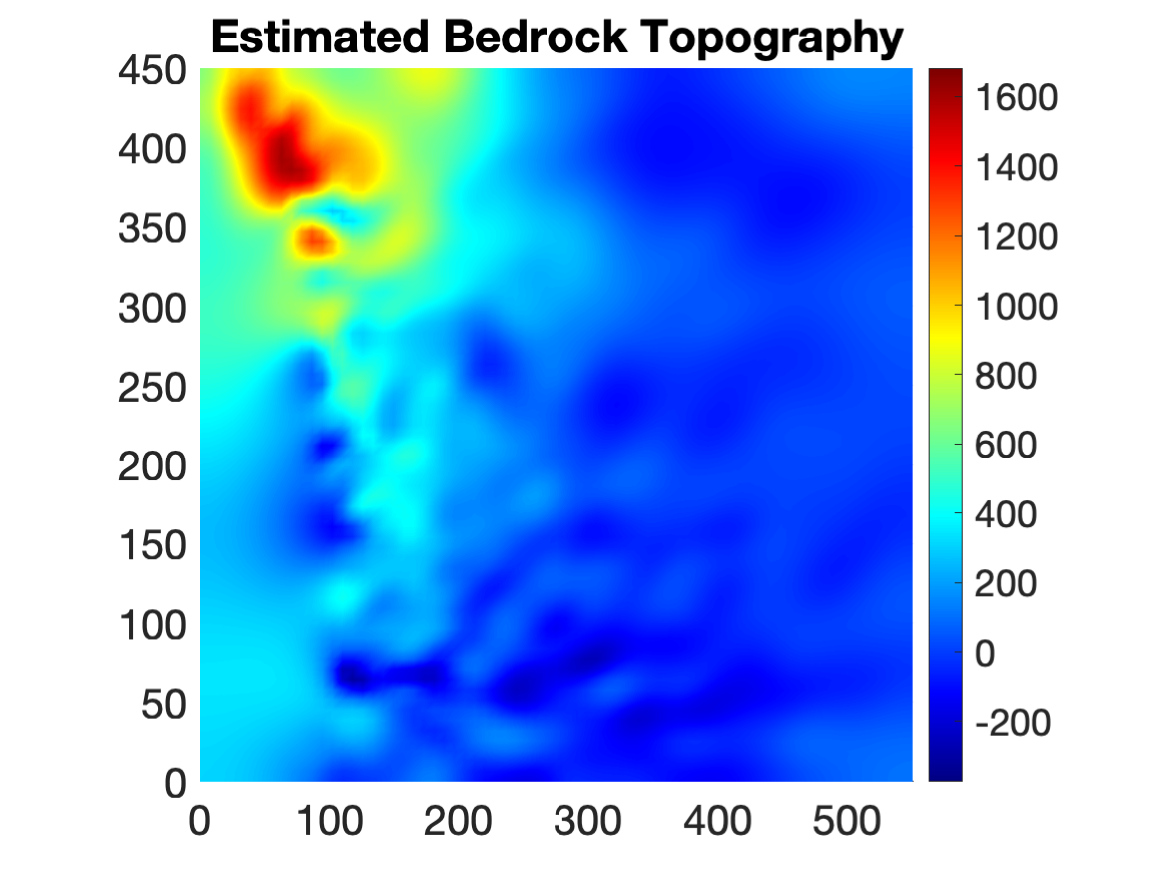}
    \caption{Left: time $t=10$ ice sheet surface height; center: bedrock topography used to synthesize data; right: bedrock topography estimated by solving the inverse problem with $\ta=\overline{\ta}$.}
  \label{fig:ice_sheet_inversion}
\end{figure}

Figure~\ref{fig:ice_sheet_inversion} displays the terminal time ($t=10$ years)
surface height, the ``true" bedrock topography used to synthesize data, and the
estimated bedrock topography that results from solving the inverse problem with
$\ta=\overline{\ta}$. The surface height is small near $x=0$ as the ice sheet is
tapering down to the coastline on the western side of Greenland. The large
bedrock topography in the top left corner corresponds to a mountain along the
coast. We observe that the inversion correctly identifies the mountain, although
it cannot resolve all of its local topography.

\subsubsection*{Preconditioned pseudo-time continuation}
The inverse problem is computationally intensive as a result of the nonlinearities of the forward model and the large data dimension. The optimization algorithm required 26 iterations using a truncated CG Newton trust region method. To better understand uncertainty due to the auxiliary parameters, we seek to solve the inverse problem for perturbations $\tilde{\ta}$. Since $\ta \in \R^{1922}$, many perturbations must be computed; consequently, solving the perturbed inverse problem solution quickly is crucial. For this article, we focus on one $\tilde{\ta}$ that perturbs the basal friction with the spatial field $\sin \left( \frac{2 \pi x}{550} \right) \sin \left( \frac{2 \pi y}{450} \right)$ and the forcing term with the field $\cos \left( \frac{2 \pi x}{550} \right) \cos \left( \frac{2 \pi y}{450} \right)$.

Given our observations from the previous example, we use the modified Euler
method as the predictor and take $N=3$ time steps. The preconditioner is
initialized using the inverse of the regularization Hessian (i.e.,
$r_\text{init}=0$) and we set $r_\text{update}=80$. This larger update rank is
determined based on our observations from the optimization iteration history
that implies the Hessian is not amenable for fast CG convergence (many
optimization iterations required $\mathcal O(100)$ CG iterations). However, we
found that only a total of $24$ vectors are stored over the three time steps as
a result of filtering (with tolerance $\tau=10^{-6}$) that reduces the number of
vectors retained for the preconditioner. 

\begin{figure}[h]
\centering
  \includegraphics[width=0.6\textwidth]{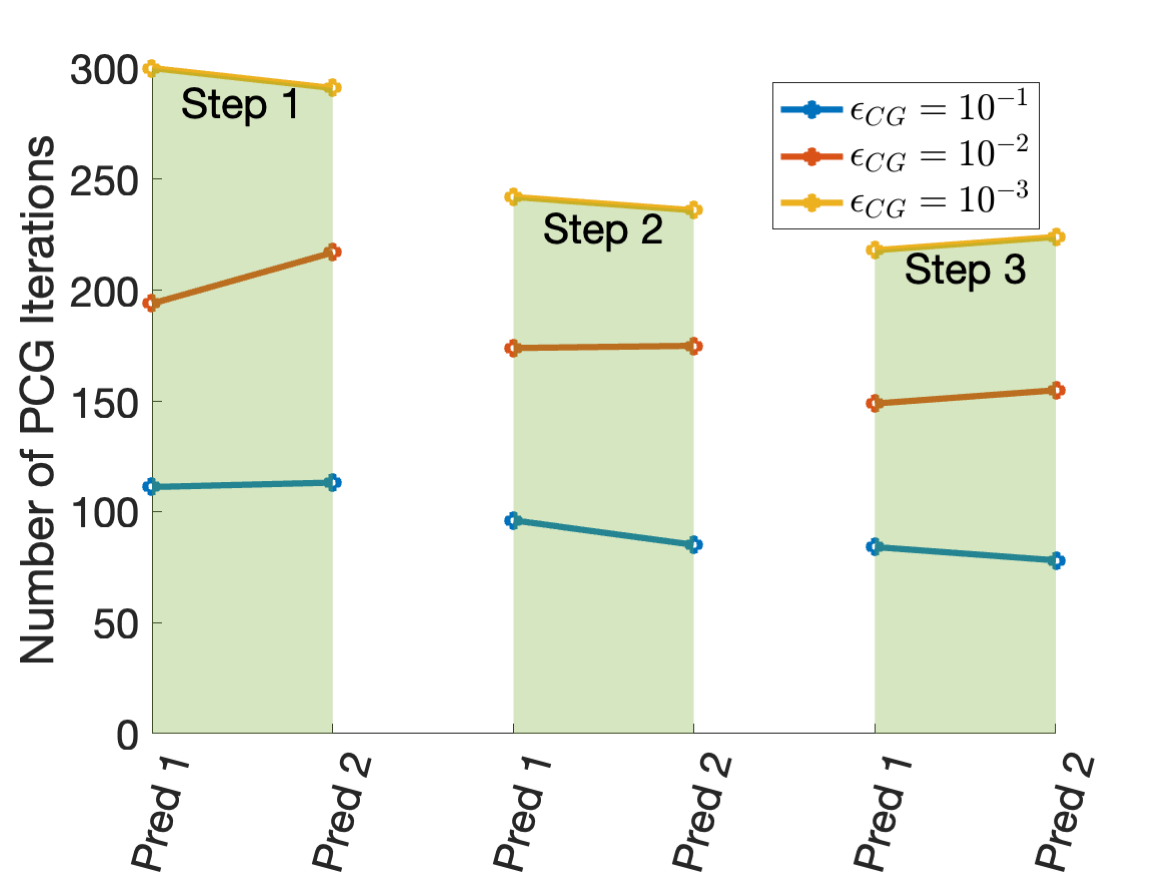}
    \caption{Summary of the preconditioned pseudo-time continuation computational cost, as measured in the number of PCG iterations. Each colored region corresponds to a time step and the points within each region indicates the number of PCG iterations required in the two predictor step Hessian inversions. A corrector step was never needed. The three lines within each step column corresponds to using the tolerances $\epsilon_\text{CG} \in \{10^{-1},10^{-2}, 10^{-3} \}$.}
  \label{fig:ice_sheet_pcg_cost}
\end{figure}

We conduct three
experiments using PCG tolerances $\epsilon_{\text{CG}}\in \{ 10^{-1}, 10^{-2},
10^{-3} \}$.  Figure~\ref{fig:ice_sheet_pcg_cost} displays the number of PCG
iterations required for the Hessian inversions. In each of the 3 time steps, we
execute two Hessian inversions in the Predictor Step. 
As noted in Section~\ref{sec:algorithm_summary}, we can check for optimality 
after the predictor step to determine whether the corrector step is necessary.
In this case, we found that no corrector steps were required.  

For a proper computational cost comparison, note that nonlinear state
solves are more expensive than the linear adjoint and incremental solves used
for gradient and Hessian-vector product computation. Each state solve was
executed using Newton's method with an average of $4$ linear solves required per
nonlinear solve. Hence, we will count cost in units of linear solves with a
state solve counting as $4$ linear solves, a gradient evaluation as $1$ linear
solve, and matrix-vector products with $\B$ and $\H$ as $2$ linear solves. Note that 
all linear solves have a comparable computational complexity since they involve
similar coefficient matrices; see Appendix~\ref{appendix:computation} for details regarding
the coefficient matrices that arise in the adjoint method.
The computational cost is summarized in
Table~\ref{tbl:computational_cost}. In this example, we observe that
the computational cost is proportional to the PCG tolerance, with a loose tolerance being superior. 

\renewcommand{\arraystretch}{1.5} 
\begin{table}[t]
    \centering
    \begin{tabular}{c|c|c|c|c|c}
 Method   &    State & Gradient & $\B$ &  $\H$ & Total Linear Solves \\
 \hline 
$\epsilon_\text{CG}=10^{-1}$ & 7 & 7 & 6 & 567 & 1181 \\
 \hline 
$\epsilon_\text{CG}=10^{-2}$ & 7 & 7 & 6 & 1064 & 2175 \\
 \hline 
$\epsilon_\text{CG}=10^{-3}$ & 7 & 7 & 6 & 1511 & 3069  \\
    \end{tabular}
    \caption{Summary of computational cost for the pseudo-time continuation algorithm using
    a modified Euler predictor, an adaptive preconditioner, and various PCG tolerances.}
    \label{tbl:computational_cost}
\end{table}

The comparison presented in Table~\ref{tbl:computational_cost} highlights that the computational cost does vary nontrivially with the selection of the PCG tolerance. However, it is noteworthy that the preconditioned pseudo-time continuation method converged with a coarse time step resolution $N=3$ and loose PCG tolerance. No algorithm parameter tuning was performed other than running with three different PCG tolerances for comparison. By contrast, we compared with re-optimization to solve the perturbed optimization problem. In this instance, the trust region radius and iteration limit for the perturbed optimization problem required tuning. We used the optimization algorithm parameters that were successful in the $\ta=\overline{\ta}$ optimization problem. However, this resulted in a large optimization step which led to a lack of convergence from the forward solve. We reduced the trust region radius and increased the maximum number of optimization iterates to ensure convergence. Such challenges in algorithm parameter specification are common in practice for nonlinear optimization problems. Since preconditioned pseudo-time continuation takes the parameter perturbation gradually, it is less prone to convergence issues that require additional algorithm parameter tuning. This robustness is particularly important in uncertainty quantification contexts when minimizers are sought for many different perturbations $\tilde{\ta}$.

\section{Conclusion}\label{sec:conc}

The preconditioned pseudo-time continuation approach is superior to executing an
optimization algorithm on the perturbed problem. The core concept is to
gradually perturb the auxiliary parameter and adapt the optimization solution
and Hessian preconditioner. This ensures desirable mathematical properties such
as the gradient equaling zero (up to a tolerance) and the Hessian being positive
definite at each time step. When combined with a predictor-corrector scheme,
this produces a computationally efficient and reliable algorithm that requires
little effort tuning algorithm parameters.

We observed that looser PCG tolerances generally yield lower computational
costs. This algorithm parameter, denoted as $\epsilon_\text{CG}$, appears to be
the most important consideration in our approach. Taking $\epsilon_\text{CG}$
too small generally increases the computational cost. Taking
$\epsilon_\text{CG}$ too large may result in a poor step and thus require
multiple Newton iterations for correction. This can increase the total
computational cost. However, the predictor-corrector method converges as long as
the time step is sufficiently small. Adaptive time stepping is another avenue of
future research to enhance the algorithms. In particular, the adaptive algorithms
of~\cite[chap. 5]{Deuflhard_book} are a promising avenue to leverage computations
within the steps to estimate Lipschitz constants needed to automatically adapt
the time step size. 

One may be tempted to use the quasi-Newton inverse Hessian approximation in
place of running PCG for the predictor step, corrector step, or both.  However,
error in the Hessian approximation may lead to running many corrector steps that
exhibit slower convergence.  Note also that in the context of optimization
problems constrained by PDEs, the linear solves required to compute
Hessian-vector products share common structures. Hence, the cost of computing
multiple Hessian-vector products can be amortized by reusing matrix
factorizations that arise from solving the incremental state and adjoint
equations.  Consequently, it is advantageous to take fewer iterations with more
Hessian-vector products per iteration rather than using an approximate Hessian
that requires many corrector iterations to ensure gradient norm tolerance satisfaction.

Our preconditioning approach reduces the computational cost measured in terms of 
the number of linearized PDE solves. However, it does introduce an additional
storage requirement that is generally not required by a PCG implementation. For
storage-limited applications, our approach may not be as effective. Nonetheless,
the filtering tolerance, denoted as $\tau$ in our article, provides a natural
mechanism to strike a balance between computational cost and storage. Furthermore, as we
demonstrated in our numerical results, a modest storage budget is often sufficient to
realize considerable benefits in computation cost reduction.

We demonstrated the benefit of using a second-order time integration method for
the predictor step. The benefit is marginal in some instances, but we observed a significant
benefit in scenarios with large auxiliary parameter
perturbations and a small number of time steps. Higher-order time integration methods may also be considered. 
Based on the insight
gained from our numerical results, we conjecture that higher-order time
integrators may provide further improvements when the auxiliary parameter
perturbation is very large. Exploration of such scenarios is an interesting line 
of future work.

Lastly, we note the potential benefit of our approach when many perturbed
optimization problems must be solved repeatedly. For instance, in applications where $N$
auxiliary parameter samples are drawn and $N$ corresponding optimization
problems must be solved. There is potential to design continuation paths and
reuse preconditioners across samples so that we may solve the $N$ perturbed
optimization problems at a computational cost that is much less than solving the
$N$ problems independently. This is a topic of ongoing research as we explore
ways to find a minimum number of paths necessary to pass through a set of $N$
auxiliary parameter samples. Using nonlinear paths may offer a considerable
reduction in the number of paths that must be followed. However, the path designs must also
consider the effects of nonlinearity in the ODE dynamics and the preconditioner updating.

\appendix

\section{Computations within a time step} \label{appendix:computation}
Within each time step, we require computation of $\vec{v}_k =
\B(\m_k,\ta_k)\Delta \ta$, $\nabla_\m J(\m_k,\ta_k)$, and $\H(\m_k,\ta_k)^{-1}
\vec{v}_k$, where the latter is computed via PCG.  In this appendix we outline
these computations using the adjoint method and discuss the
associated costs.

Assume that the parameter-to-observable map $\F(\m,\ta)$ is defined in terms
of an observation operator composed with the solution operator of a PDE.
In particular, let $\uu \in \R^m$ be the discretized state variable, $\mathcal S:\R^n \times
\R^p \to \R^m$ the discretized solution operator, and $\mathcal Q: \R^m
\to \R^d$ the observation operator. This way, $\F(\m,\ta) = \mathcal Q(\mathcal
S(\m,\ta))$. Furthermore, we define the full-space objective  
$\mathcal J:\R^m \times \R^n \to \R$ as
\begin{align*}
\mathcal J(\uu,\m) = \frac{1}{2} \| \mathcal Q(\uu) - \vec{y} \|^2 + \frac{1}{2} \m^\top \vec{R} \m
\end{align*}
so that the reduced objective $J(\m,\ta)$ in \eqref{equ:optim} satisfies $J(\m,\ta) = \mathcal J(\mathcal S(\m,\ta),\m)$.

 \begin{algorithm}[ht!!]
	\caption{Computation of $\nabla_\m J(\m_k,\ta_k)$ and $\B(\m_k,\ta_k) \Delta \ta$}
	\begin{algorithmic} [1] \label{alg:B_computation}
		\STATE \textbf{Input: } $\m_k,\ta_k$, $\Delta \ta$
		\STATE Solve the state equation $\cc(\uu,\m_k,\ta_k)=\vec{0}$ for $\uu_k=\mathcal S(\m_k,\ta_k)$
		\STATE Solve the adjoint equation $\nabla_{\uu} \cc(\uu_k,\m_k,\ta_k)^\top \vec{\lambda}_k = -\nabla_{\uu} \mathcal J(\uu_k,\m_k)$
		\STATE Compute $\nabla_\m J(\m_k,\ta_k) = \nabla_{\m} \cc(\uu_k,\m_k,\ta_k)^\top \vec{\lambda}_k + \nabla_{\m} \mathcal J(\uu_k,\m_k)$ 
		\STATE Solve the incremental state equation $\nabla_{\uu} \cc(\uu_k,\m_k,\ta_k) \vec{\xi} = -\nabla_{\ta} \cc(\uu_k,\m_k,\ta_k) \Delta \ta$
		\STATE Compute 
		$$ \vec{\nu} = \nabla_{\uu,\uu} \mathcal J(\uu_k,\m_k) \vec{\xi} + \nabla_{\uu,\uu} \cc(\uu_k,\m_k,\ta_k)[\vec{\lambda}_k] \vec{\xi} + \nabla_{\uu,\ta} \cc(\uu_k,\m_k,\ta_k)[\vec{\lambda}_k] \Delta\ta$$ 
		\STATE Solve the incremental adjoint equation $\nabla_{\uu} \cc(\uu_k,\m_k,\ta_k)^\top \vec{\beta} = -\vec{\nu}$
		\STATE Compute
		\begin{align*}
		\B(\m_k,\ta_k) \Delta \ta &= \nabla_{\m} \cc(\uu_k,\m_k,\ta_k)^\top \vec{\beta} + \nabla_{\m,\uu} \cc(\uu_k,\m_k,\ta_k)[\vec{\lambda}_k] \vec{\xi} \\
		&+  \nabla_{\m,\ta} \cc(\uu_k,\m_k,\ta_k)[\vec{\lambda}_k] \Delta \ta
		\end{align*}
	\end{algorithmic}
\end{algorithm}

Let $\cc(\uu,\m,\ta)=\vec{0} \in \R^m$ denote the discretized PDE residual such $\cc(\mathcal S(\m,\ta),\ta)$ $=\vec{0}$
that for all $(\m,\ta).$ Then we may
express the computation of  
$\nabla_\m J(\m_k,\ta_k)$ 
and 
$\vec{v}_k = \B(\m_k,\ta_k)\Delta \ta$ 
in terms of solves involving the Jacobian of $\cc$ with respect to 
$\uu$, $\nabla_{\uu} \cc$; see Algorithm~\ref{alg:B_computation}. In what follows, we denote 
the action of the PDE residual on the adjoint variable $\vec\lambda$ by 
$\cc(\uu,\m,\ta)[\vec{\lambda}] \vcentcolon= \cc(\uu,\m,\ta)^\top\vec{\lambda}$.

The computational cost of Algorithm~\ref{alg:B_computation} is approximately 1 state solve
and 3 linearized (incremental) PDE solves, since
the cost of other computations is negligible.

After executing Algorithm~\ref{alg:B_computation}, we use $\nabla_\m J(\m_k,\ta_k)$, along
with vectors computed in the previous time step, to perform the parametric
quasi-Newton update~\eqref{equ:bfgs_update}. The computation of
$\H(\m_k,\ta_k)^{-1} \B(\m_k,\ta_k) \Delta \ta$ via PCG requires $\ell$
Hessian-vector products, where
$\ell$ is the number of PCG iterations.
Algorithm~\ref{alg:H_computation} details the computation of one Hessian-vector
product, whose total cost is 2 linearized PDE solves. 

\begin{algorithm}[H]
	\caption{Computation of $\H(\m_k,\ta_k) \vec{p}$}
	\begin{algorithmic} [1] \label{alg:H_computation}
		\STATE \textbf{Input: } $\vec{p},\m_k,\ta_k,\uu_k,\vec{\lambda}_k$
		\STATE Solve the incremental state equation $\nabla_{\uu} \cc(\uu_k,\m_k,\ta_k) \vec{\xi} = -\nabla_{\m} \cc(\uu_k,\m_k,\ta_k) \vec{p}$
		\STATE Compute 
		$$ \vec{\nu} = \nabla_{\uu,\uu} \mathcal J(\uu_k,\m_k) \vec{\xi} + \nabla_{\uu,\uu} \cc(\uu_k,\m_k,\ta_k)[\vec{\lambda}_k] \vec{\xi} + \nabla_{\uu,\m} \cc(\uu_k,\m_k,\ta_k)[\vec{\lambda}_k] \vec{p}$$ 
		\STATE Solve the incremental adjoint equation $\nabla_{\uu} \cc(\uu_k,\m_k,\ta_k)^\top \vec{\beta} = -\vec{\nu}$
		\STATE Compute
		\begin{align*}
		\H(\m_k,\ta_k) \vec{p} &= \nabla_{\m} \cc(\uu_k,\m_k,\ta_k)^\top \vec{\beta} +  \nabla_{\m,\uu} \cc(\uu_k,\m_k,\ta_k)[\vec{\lambda}_k] \vec{\xi} \\
		& +  \nabla_{\m,\ta} \cc(\uu_k,\m_k,\m_k)[\vec{\lambda}_k] \vec{p} + \nabla_{\m,\m} \mathcal J(\m_k,\ta_k) \vec{p}
		\end{align*}
	\end{algorithmic}
\end{algorithm}

\section{Preconditioner initialization using low-rank approximation} \label{appendix:init_hess}
Consider the 
generalized eigenvalue problem 
\begin{equation}\label{equ:eig} 
    \HM(\tab, \ms(\tab)) \vec{v}_j = \lambda_j \HR \vec{v}_j, \quad j \in \{1, \ldots, n\}.    
\end{equation}
In ill-posed inverse problems, typically the eigenvalues decay rapidly.
Also, since we have only finite-dimensional data, the number of measurements also 
constrains the data-misfit Hessian rank. Let us define the 
weighted inner product $\ipr{\cdot}{\cdot}$ as 
\[
   \ipr{\vec{u}}{\vec{v}} = \vec{u}^\top \HR \vec{v}, \quad \vec{u}, \vec{v} \in \R^n.  
\] 
Assuming the eigenvectors are normalized, i.e., 
$\ipr{\vec{v}_i}{\vec{v}_j} = \delta_{ij}$, for $i, j \in \{1, \ldots, n\}$,
we have    
\[
\H(\m_0,\ta_0) = \HM(\m_0,\ta_0) + \HR 
 = \sum\limits_{j=1}^n \lambda_j \HR \vec{v}_j \vec{v}_j^\top \HR + \HR. 
\]
Observing that $\{ \HR^\frac{1}{2} \vec{v}_j \}_{j=1}^n$ is an orthonormal basis
and using the Sherman--Morrison--Woodbury formula, we obtain
\begin{align}
\label{equ:nominal_hessian_decom}
\H(\m_0,\ta_0)^{-1} = \HR^{-1} - \vec{V} \vec{\Gamma} \vec{V}^\top,
\end{align}
where $\vec{V} = \begin{pmatrix} \vec{v}_1 & \vec{v}_2 & \cdots & \vec{v}_n \end{pmatrix} \in \R^{n \times n}$ and $\vec{\Gamma} \in \R^{n \times n}$ is the diagonal matrix whose $(j,j)$ entry is $\lambda_j/(1+\lambda_j)$.

Given a rank $r$ truncation of~\eqref{equ:nominal_hessian_decom}, where $r$ is 
chosen such that $\lambda_r \ll \lambda_1$, we use the initial preconditioner 
\begin{align}
\label{equ:E0}
\E_0 =\HR^{-1} - \vec{V}_{\!\!r} \vec{\Gamma}_{r} \vec{V}_{\!\!r}^\top
\end{align}
where 
\[
\vec{V}_{\!\!r} =  \begin{pmatrix} \vec{v}_1 & \vec{v}_2 & \cdots & \vec{v}_r \end{pmatrix} 
\quad \text{and} \quad
\vec{\Gamma}_{r} = \mathrm{diag}\left(\frac{\lambda_1}{1+\lambda_1}, \ldots, \frac{\lambda_r}{1+\lambda_r}\right).
\]

\section*{Acknowledgments}

A.~Alexanderian was supported in part by US
National Science Foundation grant DMS \#2111044.
This work was supported by the Laboratory Directed Research and Development program at Sandia
National Laboratories, a multimission laboratory managed and operated by National Technology and
Engineering Solutions of Sandia LLC, a wholly owned subsidiary of Honeywell International Inc. for the U.S.
Department of Energy’s National Nuclear Security Administration under contract DE-NA0003525.
This written work is authored by an employee of NTESS. The employee, not NTESS, owns the right, title and interest in and to the written work and is responsible for its contents. Any subjective views or opinions that might be expressed in the written work do not necessarily represent the views of the U.S. Government. The publisher acknowledges that the U.S. Government retains a non-exclusive, paid-up, irrevocable, world-wide license to publish or reproduce the published form of this written work or allow others to do so, for U.S. Government purposes. The DOE will provide public access to results of federally sponsored research in accordance with the DOE Public Access Plan.

\bibliographystyle{siamplain}
\bibliography{refs}
\end{document}